\documentclass{amsart}
\usepackage{graphicx}
\usepackage{amscd}
\usepackage{amsmath}
\usepackage{amsfonts}
\usepackage{amssymb}
\theoremstyle{plain}
\newtheorem{theorem}{Theorem}[section]

\newtheorem{corollary}[theorem]{Corollary}
\theoremstyle{definition}
\newtheorem{definition}[theorem]{Definition}
\newtheorem{examples}[theorem]{Examples}
\theoremstyle{remark}
\newtheorem*{remark}{Remark}
\newtheorem*{remarks}{Remarks}
\numberwithin{equation}{section}

\newcommand{\EMAILHACK}{{\texttt{tonthat@math.uiowa.edu}}}
\def\openone%{\hbox{\upshape \small1\kern-3.3pt\normalsize1}}
{\mathchoice
{\hbox{\upshape \small1\kern-3.3pt\normalsize1}}
{\hbox{\upshape \small1\kern-3.3pt\normalsize1}}
{\hbox{\upshape \tiny1\kern-2.3pt\SMALL1}}
{\hbox{\upshape \Tiny1\kern-2pt\tiny1}}}
\newbox\ipbox
\newcommand{\ip}[2]{\left\langle #1\mathrel{\mathchoice
{\setbox\ipbox=\hbox{$\displaystyle \left\langle\mathstrut #1#2\right\rangle$}
\vrule height\ht\ipbox width0.25pt depth\dp\ipbox}
{\setbox\ipbox=\hbox{$\textstyle \left\langle\mathstrut #1#2\right\rangle$}
\vrule height\ht\ipbox width0.25pt depth\dp\ipbox}
{\setbox\ipbox=\hbox{$\scriptstyle \left\langle\mathstrut #1#2\right\rangle$}
\vrule height\ht\ipbox width0.25pt depth\dp\ipbox}
{\setbox\ipbox=\hbox{$\scriptscriptstyle 
\left\langle\mathstrut #1#2\right\rangle$}
\vrule height\ht\ipbox width0.25pt depth\dp\ipbox}
} #2\right\rangle}
\newcommand{\rip}[2]{\left( #1\mathrel{\mathchoice
{\setbox\ipbox=\hbox{$\displaystyle \left(\mathstrut #1#2\right)$}
\vrule height\ht\ipbox width0.25pt depth\dp\ipbox}
{\setbox\ipbox=\hbox{$\textstyle \left(\mathstrut #1#2\right)$}
\vrule height\ht\ipbox width0.25pt depth\dp\ipbox}
{\setbox\ipbox=\hbox{$\scriptstyle \left(\mathstrut #1#2\right)$}
\vrule height\ht\ipbox width0.25pt depth\dp\ipbox}
{\setbox\ipbox=\hbox{$\scriptscriptstyle \left(\mathstrut #1#2\right)$}
\vrule height\ht\ipbox width0.25pt depth\dp\ipbox}
} #2\right)}
\input cyracc.def
\font\eightcyr=wncyr8

\begin{document}
\title
[Th{\'{e}}or{\`{e}}mes de r{\'{e}}ciprocit{\'{e}}%
]{Reciprocity theorems for holomorphic representations of some 
infinite-dimensional groups\\[12pt]Quelques th{\'{e}%
}or{\`{e}}mes de r{\'{e}}ciprocit{\'{e}} pour les repr{\'{e}}%
sentations holomorphes irr{\'{e}}%
ductibles de certains groupes de dimension infinie}
\author[Tuong Ton-That]{\parbox{\textwidth}{\centering
Tuong Ton-That\\The University of Iowa\\Iowa City, IA 52242-1419 USA\\
\EMAILHACK
}}
\subjclass{02.20.Tw, 02.20.Qs, 03.65.Fd}
\begin{abstract}
Let $\mu$ denote the \emph{Gaussian measure} on $\mathbb{C}^{n\times k}$
defined by $d\mu\left(  Z\right)  =\pi^{-nk}\exp\left[  -\operatorname
*{Tr}\left(  ZZ^{\dag}\right)  \right]  \,dZ$, where $\operatorname
*{Tr}$ denotes
the trace function, $Z^{\dag}=\bar{Z}^{T}$, and $dZ$ denotes the Lebesgue
measure on $\mathbb{C}^{n\times k}$. Let $\mathcal{F}_{n\times k}$ denote the
Barg\-mann--Segal--Fock space of holomorphic entire functions on
$\mathbb{C}^{n\times k}$ which are also square-integrable with respect to
$\mu$. Fix $n$ and let $\mathcal{F}_{n\times\infty}$ denote the Hilbert-space
completion of the \emph{inductive limit} $\lim_{k\rightarrow\infty}\mathcal
{F}_{n\times k}$. Let $G_{k}$ and $H_{k}$ be compact groups such that
$H_{k}\subset G_{k}\subset\mathrm{GL}_{k}\left(  \mathbb{C}\right)  $. Let
$G_{\infty}$ (resp.\ $H_{\infty}$) denote the inductive limit $\bigcup
_{k=1}^{\infty}G_{k}$ (resp.\ $\bigcup_{k=1}^{\infty}H_{k}$). Then the
representation $R_{G_{\infty}}$ (resp.\ $R_{H_{\infty}}$) of $G_{\infty}$
(resp.\ $H_{\infty}$), obtained by right translation on
$\mathcal{F}_{n\times\infty}$, is
a \emph{holomorphic representation} of $G_{\infty}$
(resp.\ $H_{\infty}$) in the sense defined by Ol'shanskii. Then $R_{G_{\infty
}}$ and $R_{H_{\infty}}$ give rise to the dual representations
$R_{G_{n}^{\prime}}^{\prime}$ and
$R_{H_{n}^{\prime}}^{\prime}$ of the \emph{dual
pairs} $\left(  G_{n}^{\prime},G_{\infty}^{}\right)  $ and
$\left(  H_{n}^{\prime},H_{\infty}^{}\right)  $, respectively. The generalized
\emph{Bargmann--Segal--Fock space} $\mathcal{F}_{n\times\infty}$ can be
considered as both a $\left(  G_{n}^{\prime},G_{\infty}^{}\right
)  $\emph{-dual
module} and an $\left(  H_{n}^{\prime},H_{\infty}^{}\right
)  $-dual module. It is
shown that the following multiplicity-free decompositions of $\mathcal
{F}_{n\times\infty}$ into \emph{isotypic components} $\mathcal{F}%
_{n\times\infty
}=\sum\limits_{\left(  \lambda\right)  }{\!{\oplus}\,}\mathcal
{I}_{n\times\infty}^{\left(  \lambda\right)  }=
\sum\limits_{\left(  \mu\right)  }{\!{\oplus}\,}\mathcal{I}_{n\times\infty
}^{\left(  \mu\right)
}$ hold, where $\left(  \lambda\right)  $ is a common \emph{irreducible
signature} of the pair $\left(  G_{n}^{\prime},G_{\infty}^{}\right)  $ and
$\left(  \mu\right)  $ a common irreducible signature of the pair $\left(
H_{n}^{\prime},H_{\infty}^{}\right)  $, and $\mathcal{I}_{n\times\infty
}^{\left(
\lambda\right)  }$ (resp.\ $\mathcal{I}_{n\times\infty}^{\left(  \mu\right)
}$) is both the isotypic component of the equivalence classes $\left(
\lambda\right)  _{G_{\infty}}$ (resp.\ $\left(  \mu\right)  _{H_{\infty}}$)
and $\left(  \lambda^{\prime}\right)  _{G_{n}^{\prime}}$ (resp.\ $\left(
\mu^{\prime}\right)  _{H_{n}^{\prime}}$). A \emph{reciprocity theorem,} giving
the multiplicity of $\left(  \mu\right)  _{H_{\infty}}$ in the restriction to
$H_{\infty}$ of $\left(  \lambda\right)  _{G_{\infty}}$ in terms of the
multiplicity of $\left(  \lambda^{\prime}\right)  _{G_{n}^{\prime}}$ in the
restriction to $G_{n}^{\prime}$ of $\left(  \mu^{\prime}\right)
_{H_{n}^{\prime}}$, constitutes the main result of this paper. Several
applications of this theorem to Physics are also discussed.\\[10pt]
\textsc{R{\'{e}}sum{\'{e}}.} Soit $\mu
$ la mesure de Gauss definie sur l'espace vectoriel $\mathbb{C}^{n\times k}$
par la formule
\newline
\begin{minipage}{\textwidth}
\hskip-36pt\begin{minipage}{\textwidth}
\begin{equation*}
d\mu\left(  Z\right)  =\pi^{-nk}\exp\left[  -\operatorname
*{Tr}\left(  ZZ^{\dag}\right)  \right]  \,dZ,\qquad z\in\mathbb{C}^{n\times
k},
\end{equation*}
\end{minipage}\hskip-36pt\vspace*{\abovedisplayskip}
\end{minipage}
o{\`{u}} l'on d{\'{e}}signe par $\operatorname
*{Tr}$
la trace d'une matrice, $Z^{\dag}=\bar{Z}^{T}%
$, et par $dZ$ la mesure de Lebesgue
sur $\mathbb{C}^{n\times k}$. Soit $\mathcal{F}_{n\times k}%
$ l'espace hilbertien de
Barg\-mann--Segal--Fock des fonctions enti{\`{e}}res holomorphes
$f\colon\mathbb{C}^{n\times k}\rightarrow\mathbb{C}$ telles que
$f$ soient de carr{\'{e}}-integrable par rapport {\`{a}}
la mesure $\mu$. On fixe $n$ et l'on d{\'{e}}signe par
$\mathcal{F}_{n\times\infty}$ le
compl{\'{e}}t{\'{e}} de la \emph{limite inductive} par
rapport {\`{a}} $k$ des espaces $\mathcal
{F}_{n\times k}$. Pour chaque
$k$ soient $G_{k}$ et $H_{k}$ deux groupes compacts tels que
$H_{k}\subset G_{k}\subset\mathrm{GL}_{k}\left(  \mathbb{C}\right
)  $, et l'on suppose aussi que
$H_{k-1}\subset H_{k}\subset H_{k+1}\subset\cdots$ et
$G_{k-1}\subset G_{k}\subset G_{k+1}\subset\cdots$. Soit
$G_{\infty}$ (resp.\ $H_{\infty}$)  la limite inductive de la chaine
$\left\{
G_{k}\right\}$ (resp.\ $\left\{ H_{k}\right\}$). Alors la
repr{\'{e}}sentation $R_{G_{\infty}}$ (resp.\ $R_{H_{\infty}}$) de $G_{\infty
}$
(resp.\ $H_{\infty}$),
obtenue par translation {\`{a}} droite sur
$\mathcal{F}_{n\times\infty}$, est \emph{holomorphe}
dans le sens de Ol'shanskii. Les repr{\'{e}}sentations $R_{G_{\infty
}}$ et $R_{H_{\infty}}$ donnent
lieu aux repr{\'{e}}sentations
$R_{G_{n}^{\prime}}^{\prime}$ et
$R_{H_{n}^{\prime}}^{\prime}$, respectivement, des \emph{paires
duales}
$\left(  G_{n}^{\prime},G_{\infty}^{}\right)  $ et
$\left(  H_{n}^{\prime},H_{\infty}^{}\right)  $. L'espace
hilbertien generalis{\'{e}}
de Bargmann--Segal--Fock $\mathcal{F}_{n\times\infty}$ peut {\^{e}}tre
consider{\'{e}} en m{\^{e}}me temps comme un
$\left(  G_{n}^{\prime},G_{\infty}^{}\right)  $-module et un
$\left(  H_{n}^{\prime},H_{\infty}^{}\right)  $-module. On
montre que l'on a les d{\'{e}}compositions suivantes de $\mathcal
{F}_{n\times\infty}$ en uniques \emph{composantes isotypiques}
\newline
\begin{minipage}{\textwidth}
\hskip-36pt\begin{minipage}{\textwidth}
\begin{equation*}
\mathcal{F}_{n\times\infty
}=\sum\limits_{\left(  \lambda\right)  }{\!{\oplus}\,}\mathcal
{I}_{n\times\infty}^{\left(  \lambda\right)  }=
\sum\limits_{\left(  \mu\right)  }{\!{\oplus}\,}
\mathcal{I}_{n\times\infty}^{\left(  \mu\right)  },
\end{equation*}
\end{minipage}\hskip-36pt\vspace*{\abovedisplayskip}
\end{minipage}
o{\`{u}} $\left(  \lambda\right)  $ est une \emph{signature irr{\'{e}%
}ductible}
commune de la paire $\left(  G_{n}^{\prime},G_{\infty}^{}\right)  $ et
$\left(  \mu\right)  $ celle de la paire $\left(
H_{n}^{\prime},H_{\infty}^{}\right)  $, et o{\`{u}} $\mathcal{I}%
_{n\times\infty}^{\left(
\lambda\right)  }$ (resp.\ $\mathcal{I}_{n\times\infty}^{\left(  \mu\right)
}$) est {\`{a}} la fois la composante isotypique de la
classe d'{\'{e}}quivalence de $\left(
\lambda\right)  _{G_{\infty}}$ (resp.\ $\left(  \mu\right)  _{H_{\infty}}$)
et celle de $\left(  \lambda^{\prime}\right)  _{G_{n}^{\prime}}%
$ (resp.\ $\left(
\mu^{\prime}\right)  _{H_{n}^{\prime}}$). On donne une d{\'{e}}%
monstration d'un
\emph{th{\'{e}}or{\`{e}}me de r{\'{e}}ciprocit{\'{e}},} donnant
la multiplicit{\'{e}} de $\left(  \mu\right)  _{H_{\infty}}%
$ dans la restriction {\`{a}}
$H_{\infty}$ de $\left(  \lambda\right)  _{G_{\infty}}$, en fonction de la
multiplicit{\'{e}} de $\left(  \lambda^{\prime}\right)  _{G_{n}^{\prime}%
}$ dans la
restriction {\`{a}} $G_{n}^{\prime}$ de $\left(  \mu^{\prime}\right)
_{H_{n}^{\prime}}$. L'article se termine par une discussion de plusieurs
applications en Physique du th{\'{e}}or{\`{e}}me pr{\'{e}}c{\'{e}}dant.
\end{abstract}
\maketitle
%EndExpansion

\section{\label{Int}Introduction}

In recent years there is great interest, both in Physics and in Mathematics,
in the theory of unitary representations of infinite-dimensional groups and
their Lie algebras (see, for example, \cite{Kac90}, and the literature cited
therein). Starting with the seminal work of I.~Segal in \cite{Seg50} the
representation theory of $\mathrm{U}\left(  \infty\right)  $ and other
classical infinite-dimensional groups was thoroughly investigated by Kirillov
in \cite{Kir73}, Stratila and Voiculescu in \cite{StVo75}, Pickrell in
\cite{Pic87}, Ol'shanskii in \cite{Ol'88}, Gelfand and Graev in \cite{GeGr90},
Kac in \cite{Kac80}, to cite just a few. A more complete list of references
can be found in the comprehensive and important work of Ol'shanskii in
\cite{Ol'90}.

In \cite{Ol'90} Ol'shanski generalized Howe's theory of dual pairs to some
infinite-dimensional dual pairs of groups. Recently in \cite{Ton98} and
\cite{Ton97} we investigated the generalized Casimir invariants of these
infinite-dimensional dual pairs. In \cite{Ton95} we gave a general reciprocity
theorem for finite-dimensional dual pairs of groups which generalized our
previous results in \cite{KlTo89} and \cite{LeTo95}. In this article we give a
generalization of this reciprocity theorem to the case of dual pairs where one
member is infinite-dimensional and the other is finite-dimensional, and
discuss the general case where both members are infinite-dimensional. If
Section \ref{FD} we will review the reciprocity theorem given in \cite{Ton95}
which serves as the necessary background for the generalized theorem, and more
importantly, discuss several interesting applications of this theorem. Section
\ref{FID} deals with our main theorem, and the paper ends with a short
conclusion in Section \ref{Con}.

\section{\label{FD}The Reciprocity Theorem for Finite-Dimensional Pairs of
Groups and Its Applications}

In \cite{Ton95} our reciprocity theorem can be applied to the more general
context of dual representations but for this paper we shall restrict ourself
to the case of the oscillator dual representations and where one of the
members is a compact group.

Let $\mathbb{C}^{n\times k}$ denote the vector space of all $n\times k$
complex matrices. Let $\mu$ denote the Gaussian measure on $\mathbb{C}%
^{n\times k}$ defined by%
\begin{equation}
d\mu\left(  Z\right)  =\pi^{-nk}\exp\left[  -\operatorname*{Tr}\left(
ZZ^{\dag}\right)  \right]  \,dZ,\qquad Z\in\mathbb{C}^{n\times k},
\label{eqFD.1}%
\end{equation}
where in Eq.\ (\ref{eqFD.1}) $Z^{\dag}$ denotes the adjoint of the matrix $Z$
and $dZ$ denotes the Lebesgue measure on $\mathbb{C}^{n\times N}$. Let
$\mathcal{F}_{n\times k}\equiv\mathcal{F}\left(  \mathbb{C}^{n\times
k}\right)  $ denote the Bargmann--Segal--Fock space of all holomorphic entire
functions on $\mathbb{C}^{n\times k}$ which are also square-integrable with
respect to $d\mu$. Endowed with the inner product%
\begin{equation}%
%TCIMACRO{\TeXButton{rip}{\rip{f}{g}}}%
%BeginExpansion
\rip{f}{g}%
%EndExpansion
=\int_{\mathbb{C}^{n\times k}}f\left(  Z\right)  \overline{g\left(  Z\right)
}\,d\mu\left(  Z\right)  \mathrel{;}\qquad f,g\in\mathcal{F}_{n\times k},
\label{eqFD.2}%
\end{equation}
$\mathcal{F}_{n\times k}$ has a Hilbert-space structure. It can be easily
verified that the inner product%
\begin{equation}
\ip{f}{g}=f\left(  D\right)  \overline{g\left(  \bar{Z}\right)  }|_{Z=0}
\label{eqFD.3}%
\end{equation}
where $f\left(  D\right)  $ denotes the formal power series obtained by
replacing $Z_{\alpha j}$ by the partial derivative $\partial/\partial
Z_{\alpha j}$ ($1\leq\alpha\leq n$, $1\leq j\leq k$). In fact if $\left(
r\right)  =\left(  r_{11},\dots,r_{nk}\right)  $ is a multi-index of integers
$r_{\alpha j}\geq0$ let $Z^{\left(  r\right)  }\equiv Z_{11}^{r_{11}}\cdots
Z_{nk}^{r_{nk}}$ and $\left(  r\right)  !=r_{11}!\cdots r_{nk}!$ then it is
easy to verify that%
\begin{equation}%
%TCIMACRO{\TeXButton{rip}{\rip{\frac{Z^{\left( r\right) }}{\left[ \left
%( r\right) !\right] ^{\frac{1}{2}}}}{\frac{Z^{\left( r^{\prime}\right) }%
%}{\left[ \left( r^{\prime}\right) !\right] ^{\frac{1}{2}}}}}}%
%BeginExpansion
\rip{\frac{Z^{\left( r\right) }}{\left[ \left
( r\right) !\right] ^{\frac{1}{2}}}}{\frac{Z^{\left( r^{\prime}\right) }%
}{\left[ \left( r^{\prime}\right) !\right] ^{\frac{1}{2}}}}%
%EndExpansion
=\ip{\frac{Z^{\left( r\right) }}{\left[ \left( r\right) !\right] ^{\frac{1}%
{2}}}}{\frac{Z^{\left( r^{\prime}\right) }}{\left[ \left( r^{\prime}%
\right) !\right] ^{\frac{1}{2}}}}=\delta_{\left(  r\right)  \left(  r^{\prime
}\right)  }. \label{eqFD.4}%
\end{equation}
It follows immediately from Eq.\ (\ref{eqFD.4}) that $\left\{  Z^{\left(
r\right)  }\mathop{\big/}\left[  \left(  r\right)  !\right]  ^{\frac{1}{2}%
}\right\}  _{\left(  r\right)  }$ forms an orthonormal basis for
$\mathcal{F}_{n\times k}$ when $\left(  r\right)  $ ranges over all
multi-indices; moreover $\mathcal{P}_{n\times k}\equiv\mathcal{P}\left(
\mathbb{C}^{n\times k}\right)  $, the subspace of all polynomial functions on
$\mathbb{C}^{n\times k}$, is dense in $\mathcal{F}_{n\times k}$.

Let $G$ and $G^{\prime}$ be two topological groups. Let $R_{G}$ and
$R_{G^{\prime}}^{\prime}$ be continuous unitary and completely (discretely)
reducible representations of $G$ and $G^{\prime}$ on $\mathcal{F}_{n\times k}$
such that $R_{G}$ and $R_{G^{\prime}}^{\prime}$ commute. Then we have the
following definition of \emph{dual representations} (for the definition of
dual representations in a more general context see \cite{Ton95}).

\begin{definition}
\label{DefFD.1}The representations $R_{G}$ and $R_{G^{\prime}}^{\prime}$ are
said to be \emph{dual} if the $G^{\prime}\times G$-module $\mathcal{F}%
_{n\times k}$ is decomposed into a multiplicity-free orthogonal direct sum of
the form%
\begin{equation}
\mathcal{F}_{n\times k}=\sum\limits_{\left(  \lambda\right)  }{\!{\oplus}%
\,}\mathcal{I}_{n\times k}^{\left(  \lambda\right)  }, \label{eqFD.5}%
\end{equation}
where in Eq.\ \textup{(\ref{eqFD.5})} the label $\left(  \lambda\right)  $
characterizes both an equivalence class of an irreducible unitary
representation $\lambda_{G}$ of $G$ and an equivalence class of an irreducible
representation $\lambda_{G^{\prime}}^{\prime}$, and $\mathcal{I}_{n\times
k}^{\left(  \lambda\right)  }$ denotes the $\left(  \lambda\right)
$\emph{-isotypic component,} i.e., the direct sum \textup{(}not
canonical\/\textup{)} of all irreducible subrepresentations of $R_{G}$
\textup{(}resp.\ $R_{G^{\prime}}^{\prime}$\textup{)} that belong to the
equivalence class $\lambda_{G}$ \textup{(}resp.\ $\lambda_{G^{\prime}}%
^{\prime}$\textup{).} Moreover the $G^{\prime}\times G$-submodule
$\mathcal{I}_{n\times k}^{\left(  \lambda\right)  }$ is irreducible for all
\emph{signatures} $\left(  \lambda\right)  $; i.e., $\mathcal{I}_{n\times
k}^{\left(  \lambda\right)  }\approx V^{\left(  \lambda_{G}\right)  }%
\mathop{\hat{\otimes}}W^{\left(  \lambda_{G^{\prime}}^{\prime}\right)  }$,
where $V^{\left(  \lambda_{G}\right)  }$ \textup{(}resp.\ $W^{\left(
\lambda_{G^{\prime}}^{\prime}\right)  }$\textup{)} is an irreducible
$G$-module of class $\left(  \lambda_{G}\right)  $ \textup{(}resp.\ $G^{\prime
}$-module of class $\left(  \lambda_{G^{\prime}}^{\prime}\right)  $\textup{).}

We refer to the decomposition \textup{(\ref{eqFD.5})} as the \emph{canonical
decomposition} of the $G^{\prime}\times G$-module $\mathcal{F}_{n\times k}$.
\end{definition}

In this context we have the following theorem which is a special case of
Theorem \textup{4.1} in \cite{Ton95}.

\begin{theorem}
\label{ThmFD.2}Let $G$ be a compact group. Let $R_{G}$ and $R_{G^{\prime}%
}^{\prime}$ be given dual representations on $\mathcal{F}_{n\times k}$. Let
$H$ be a compact subgroup of $G$ and let $R_{H}$ be the representation of $H$
on $\mathcal{F}_{n\times k}$ obtained by restricting $R_{G}$ to $H$. If there
exists a group $H^{\prime}\supset G^{\prime}$ and a representation
$R_{H^{\prime}}^{\prime}$ on $\mathcal{F}_{n\times k}$ such that
$R_{H^{\prime}}^{\prime}$ is dual to $R_{H}$ and $R_{G^{\prime}}^{\prime}$ is
the restriction of $R_{H^{\prime}}^{\prime}$ to the subgroup $G^{\prime}$ of
$H^{\prime}$ then we have the following multiplicity-free decompositions of
$\mathcal{F}_{n\times k}$ into isotypic components%
\begin{equation}
\mathcal{F}_{n\times k}=\sum\limits_{\left(  \lambda\right)  }{\!{\oplus}%
\,}\mathcal{I}_{n\times k}^{\left(  \lambda\right)  }= \sum\limits_{\left(
\mu\right)  }{\!{\oplus}\,}\mathcal{I}_{n\times k}^{\left(  \mu\right)  }
\label{eqFD.6}%
\end{equation}
where $\left(  \lambda\right)  $ is a common irreducible signature of the pair
$\left(  G^{\prime},G\right)  $ and $\left(  \mu\right)  $ is a common
irreducible signature of the pair $\left(  H^{\prime},H\right)  $.

If $\lambda_{G}$ \textup{(}resp.\ $\lambda_{G^{\prime}}^{\prime}$\textup{)}
denotes an irreducible unitary representation of class $\left(  \lambda
\right)  $ \linebreak and $\mu^{}_{H}$ \textup{(}resp.\ $\mu_{H^{\prime}%
}^{\prime}$\textup{)} denotes an irreducible unitary representation of class
$\left(  \mu\right)  $ \linebreak then the multiplicity $\dim\left[
\operatorname*{Hom}_{H}\left(  \mu^{}_{H}:\lambda_{G}|_{H}\right)  \right]  $
of the irreducible representation \linebreak $\mu^{}_{H}$ in the restriction
to $H$ of the representation $\lambda_{G}$ is equal to the multiplicity
$\dim\left[  \operatorname*{Hom}_{G^{\prime}}\left(  \lambda_{G^{\prime}%
}^{\prime}:\mu_{H^{\prime}}^{\prime}|_{G^{\prime}}\right)  \right]  $ of the
irreducible representation $\lambda_{G^{\prime}}^{\prime}$ in the restriction
to $G^{\prime}$ of the representation $\mu_{H^{\prime}}^{\prime}$.
\end{theorem}

\begin{remarks}
In many cases $\operatorname*{Hom}_{H}\left(  \mu:\lambda_{G}|_{H}\right)  $
and $\operatorname*{Hom}_{G}\left(  \lambda_{G}^{\prime}:\mu_{H^{\prime}%
}^{\prime}|_{G}\right)  $ are shown to be isomorphic and can be explicitly
constructed in terms of generalized Casimir operators as given in
\cite{KlTo92} and \cite{LeTo94}.
\end{remarks}

To illustrate this theorem we devote the rest of this section to some typical
examples and discuss their generalization.

\begin{examples}
\label{ExaFD.3} \quad
%\begin{enumerate}
%\item \label{ExaFD.3(1)} 
%
%
%
%
%
%
1) Consider $\mathcal{F}_{1\times k}$ with $k\geq2$; then $\mathcal{F}%
_{1\times k}$ is the classical Barg\-mann space first considered by
V.~Bargmann in \cite{Bar61}. Then $\mathcal{P}_{1\times k}$ is the algebra of
all polynomial functions in $k$ variables $\left(  Z_{1},\dots,Z_{k}\right)
=Z$. Let $G=\mathrm{U}\left(  k\right)  $ and $G^{\prime}=\mathrm{U}\left(
1\right)  $; then the complexification of $\mathrm{U}\left(  k\right)  $
\textup{(}resp.\ $\mathrm{U}\left(  1\right)  $\textup{)} is $G_{\mathbb{C}%
}=\mathrm{GL}_{k}\left(  \mathbb{C}\right)  $ \textup{(}resp.\ $G_{\mathbb{C}%
}^{\prime}=\mathrm{GL}_{1}\left(  \mathbb{C}\right)  $\textup{).} An element
$f$ of $\mathcal{F}_{1\times k}$ is of the form%
\begin{equation}
f\left(  Z\right)  =\sum_{\left|  \left(  r\right)  \right|  =0}^{\infty
}c_{\left(  r\right)  }Z^{\left(  r\right)  } \label{eqFD.7}%
\end{equation}
with $\left(  r\right)  =\left(  r_{1},\dots,r_{k}\right)  $, $\left|  \left(
r\right)  \right|  =r_{1}+\dots+r_{k}$, and $Z^{\left(  r\right)  }%
=Z_{1}^{r_{1}}\cdots Z_{k}^{r_{k}}$, $c_{\left(  r\right)  }\in\mathbb{C}$
such that $\sum_{\left|  \left(  r\right)  \right|  =0}^{\infty}\left|
c_{\left(  r\right)  }\right|  ^{2}\left(  r\right)  !<\infty$, where $\left(
r\right)  !=r_{1}!\cdots r_{k}!$. The system $\left\{  Z^{\left(  r\right)
}\mathop{\big/}\left[  \left(  r\right)  !\right]  ^{\frac{1}{2}}\right\}  $,
where $\left(  r\right)  $ ranges over all multi-indices, forms an orthonormal
basis for $\mathcal{F}_{1\times k}$. $R_{G_{\mathbb{C}}}$ and $R_{G}$ are
defined by%
\begin{equation}%
%TCIMACRO{\TeXButton{cases}{\begin{cases}
%\left[ R_{G_{\mathbb{C}}}\left( g\right) f\right] \left( Z\right
%) =f\left( Zg\right) , & g\in\mathrm{GL}_{k}\left( \mathbb{C}\right) ,  \\
%\left[ R_{G}\left( u\right) f\right] \left( Z\right) =f\left( Zu\right
%) , & u\in\mathrm{U}\left( k\right) .
%\end{cases}}}%
%BeginExpansion
\begin{cases}
\left[ R_{G_{\mathbb{C}}}\left( g\right) f\right] \left( Z\right
) =f\left( Zg\right) , & g\in\mathrm{GL}_{k}\left( \mathbb{C}\right) ,  \\
\left[ R_{G}\left( u\right) f\right] \left( Z\right) =f\left( Zu\right
) , & u\in\mathrm{U}\left( k\right) .
\end{cases}%
%EndExpansion
\label{eqFD.8}%
\end{equation}
$R_{G_{\mathbb{C}}^{\prime}}^{\prime}$ and $R_{G^{\prime}}^{\prime}$ are
defined by%
\begin{equation}%
%TCIMACRO{\TeXButton{cases}{\begin{cases}
%\left[ R^{\prime}_{G^{\prime}_{\mathbb{C}}}\left( g^{\prime}\right
%) f\right] \left( Z\right) =f\left( \left( g^{\prime}\right) ^{t}%
%Z\right) , & g^{\prime}\in\mathrm{GL}_{1}\left( \mathbb{C}\right) ,  \\
%\left[ R^{\prime}_{G^{\prime}}\left( u^{\prime}\right) f\right] \left( Z\right
%) =f\left( \left( u^{\prime}\right) ^{t}Z\right) , & u^{\prime}\in\mathrm
%{U}\left( 1\right) .
%\end{cases}}}%
%BeginExpansion
\begin{cases}
\left[ R^{\prime}_{G^{\prime}_{\mathbb{C}}}\left( g^{\prime}\right
) f\right] \left( Z\right) =f\left( \left( g^{\prime}\right) ^{t}%
Z\right) , & g^{\prime}\in\mathrm{GL}_{1}\left( \mathbb{C}\right) ,  \\
\left[ R^{\prime}_{G^{\prime}}\left( u^{\prime}\right) f\right] \left( Z\right
) =f\left( \left( u^{\prime}\right) ^{t}Z\right) , & u^{\prime}\in\mathrm
{U}\left( 1\right) .
\end{cases}%
%EndExpansion
\label{eqFD.9}%
\end{equation}

The infinitesimal action of $R_{G_{\mathbb{C}}}$ is given by%
\begin{equation}
R_{ij}=Z_{i}\frac{\partial\;}{\partial Z_{j}},\qquad1\leq i,j\leq k,
\label{eqFD.10}%
\end{equation}
which form a basis for a Lie algebra isomorphic to \textrm{gl}$_{k}\left(
\mathbb{C}\right)  $.

The infinitesimal action of $R_{G_{\mathbb{C}}^{\prime}}^{\prime}$ is given by%
\begin{equation}
L=\sum_{i=1}^{k}Z_{i}\frac{\partial\;}{\partial Z_{i}}, \label{eqFD.11}%
\end{equation}
which forms a basis for a Lie algebra isomorphic to $\mathrm{gl}_{1}\left(
\mathbb{C}\right)  $. If $p,q\in\mathcal{P}_{1\times k}$ then from
Eq.\ \textup{(2.1)} of \cite{Ton76a} we have%
\begin{equation}
R_{G_{\mathbb{C}}}\left(  g\right)  p\left(  D\right)  R_{G_{\mathbb{C}}%
}\left(  g^{-1}\right)  =\left[  R_{G_{\mathbb{C}}}\left(  g^{\checkmark
}\right)  p\right]  \left(  D\right)  ,\qquad g\in\mathrm{GL}_{k}\left(
\mathbb{C}\right)  ,\;g^{\checkmark}=\left(  g^{-1}\right)  ^{t},
\label{eqFD.12}%
\end{equation}
so that if $u\in\mathrm{U}\left(  k\right)  $ then%
\begin{align}
\ip{R_{G}\left( u\right) p}{R_{G}\left( u\right) q}  &  =\left[  R_{G}\left(
u\right)  p\right]  \left(  D\right)  \overline{\left(  R_{G}\left(  u\right)
q\right)  \left(  \bar{Z}\right)  }\bigg|_{Z=0}\label{eqFD.13}\\
&  =R_{G}\left(  u^{\checkmark}\right)  p\left(  D\right)  R_{G}\left(
u^{t}\right)  R\left(  \bar{u}\right)  \overline{q\left(  \bar{Z}\right)
}\bigg|_{Z=0}\nonumber\\
&  =p\left(  D\right)  R_{G}\left(  u^{t}\bar{u}\right)  q\left(  \bar
{Z}u^{\checkmark}\right)  \bigg|_{Z=0}\nonumber\\
&  =\ip{p}{q},\nonumber
\end{align}
since $u^{t}\bar{u}=1$. A similar computation shows that%
\begin{equation}
R_{G_{\mathbb{C}}^{\prime}}^{\prime}\left(  g^{\prime}\right)  p\left(
D\right)  R_{G_{\mathbb{C}}^{\prime}}^{\prime}\left(  \left(  g^{\prime
}\right)  ^{-1}\right)  =\left[  R\left(  \left(  g^{\prime}\right)
^{\checkmark}\right)  \right]  \left(  D\right)  ,\qquad g^{\prime}%
\in\mathrm{GL}_{1}\left(  \mathbb{C}\right)  , \label{eqFD.14}%
\end{equation}
so that if $u\in\mathrm{U}\left(  1\right)  $ then%
\begin{equation}
\ip{R^{\prime}_{G^{\prime}}\left( u^{\prime}\right) p}{R^{\prime}_{G^{\prime}%
}\left( u^{\prime}\right) q}=\ip{p}{q}. \label{eqFD.15}%
\end{equation}
Note that all equations above from \textup{(\ref{eqFD.12})} to
\textup{(\ref{eqFD.15})} remain valid if we replace $\mathbb{C}^{1\times k}$
by $\mathbb{C}^{n\times k}$ and $\mathrm{GL}_{1}\left(  \mathbb{C}\right)  $
\textup{(}resp.\ $\mathrm{U}\left(  1\right)  $\textup{)} by $\mathrm{GL}%
_{n}\left(  \mathbb{C}\right)  $ \textup{(}resp.\ $\mathrm{U}\left(  n\right)
$\textup{).}

It follows that $R_{G}$, $G=\mathrm{U}\left(  k\right)  $ \textup{(}%
resp.\ $R_{G^{\prime}}^{\prime}$, $G^{\prime}=\mathrm{U}\left(  n\right)
$\textup{)} is a continuous unitary representation of $G$ \textup{(}%
resp.\ $G^{\prime}$\textup{)} on $\mathcal{F}_{n\times k}$.

Let $\mathcal{P}_{1\times k}^{\left(  m\right)  }$ denote the subspace
\textup{(}of $\mathcal{F}_{1\times k}$\textup{)} of all homogeneous polynomial
functions of degree $m\geq0$. Then by the Borel--Weil theorem \textup{(}see,
e.g., \cite{Ton76a}\textup{)} the restriction of $R_{G_{\mathbb{C}}}$ to
$\mathcal{P}_{1\times k}^{\left(  m\right)  }$ is an irreducible
subrepresentation of $R_{G_{\mathbb{C}}}$ with highest weight $(\underset
{k}{\underbrace{m,0,\dots,0}})$ and highest weight vector $cZ_{1}^{m}$,
$c\in\mathbb{C}^{\ast}$. In fact, by letting the infinitesimal operators
$R_{ij}$ act on $\mathcal{P}_{1\times k}^{\left(  m\right)  }$ one can easily
show that $\mathcal{P}_{1\times k}^{\left(  m\right)  }$ is an irreducible
subrepresentation of $R_{G_{\mathbb{C}}}$. By ``Weyl's unitarian trick'' the
restriction of this irreducible subrepresentation to $G$ gives an irreducible
unitary representation of $G$.

Let $0\neq p\in\mathcal{P}_{1\times k}^{\left(  m\right)  }$. Then $\left(
R_{G^{\prime}}^{\prime}\left(  g^{\prime}\right)  p\right)  \left(  Z\right)
=p\left(  \left(  g^{\prime}\right)  ^{t}Z\right)  =p\left(  g^{\prime
}Z\right)  =\left(  g^{\prime}\right)  ^{m}p\left(  Z\right)  $ for all
$g^{\prime}\in\mathrm{GL}_{1}\left(  \mathbb{C}\right)  $. So the
one-dimensional subspace of $\mathcal{F}_{1\times k}$ spanned by $p$ is an
irreducible $G_{\mathbb{C}}^{\prime}$-submodule with highest weight $\left(
m\right)  $ and its restriction to $G^{\prime}$ is an irreducible unitary
$G^{\prime}$-submodule. In fact, Euler's formula implies that%
\begin{equation}
Lp=mp\text{,\qquad for all }p\in\mathcal{P}_{1\times k}^{\left(  m\right)  }.
\label{eqFD.16}%
\end{equation}
Thus the canonical decomposition of the $G^{\prime}\times G$-module
$\mathcal{F}_{1\times k}$ is simply
\begin{equation}
\mathcal{F}_{1\times k}=\sum\limits_{\makebox[0pt]{\hss$\scriptstyle m=0$\hss
}}^{\infty}{\!{\oplus}\,}\mathcal{P}_{1\times k}^{\left(  m\right)  }.
\label{eqFD.17}%
\end{equation}
Let $H$ denote the special orthogonal subgroup $\mathrm{SO}\left(  k\right)
$. Then $H_{\mathbb{C}}=\mathrm{SO}_{k}\left(  \mathbb{C}\right)  $. Then the
ring of all $H$ \textup{(}or $H_{\mathbb{C}}$\textup{)}-invariant polynomials
in $\mathcal{P}_{1\times k}$ is generated by the constants and $p_{0}\left(
Z\right)  =\sum_{1\leq i\leq k}Z_{i}^{2}$. The ring of all $H$ \textup{(}or
$H_{\mathbb{C}}$\textup{)}-invariant differential operators with constant
coefficients is generated by the constants and the Laplacian $\bigtriangleup
=p_{0}\left(  D\right)  =\sum_{1\leq i\leq k}\partial^{2}/\partial Z_{i}^{2}$.
To find the dual representation of $R_{H}$ we follow the method given in
\cite{Ton95} by setting%
\begin{equation}
X^{+}=\frac{1}{2}p_{0},\qquad X^{-}=\frac{1}{2}p_{0}\left(  D\right)
=\frac{1}{2}\bigtriangleup\text{,\qquad and }E=\frac{k}{2}+L. \label{eqFD.18}%
\end{equation}
Then $X^{+}$ \textup{(}resp.\ $X^{-}$\textup{)} acts on $\mathcal{F}_{1\times
k}$ as a \emph{creation} \textup{(}resp.\ \emph{annihilation}\textup{)}
\emph{operator} and $E$ acts on $\mathcal{F}_{1\times k}$ as a number
operator. In fact, if $p\in\mathcal{P}_{1\times k}^{\left(  m\right)  }$ then
$X^{+}p=\frac{1}{2}p_{0}p$, $X^{-}p=\frac{1}{2}\bigtriangleup p$, and
$Ep=\left(  \left(  k/2\right)  +m\right)  p$, so that $X^{+}$ \emph{raises}
$\mathcal{P}_{1\times k}^{\left(  m\right)  }$ to $\mathcal{P}_{1\times
k}^{\left(  m+2\right)  }$, $X^{-}$ \emph{lowers} $\mathcal{P}_{1\times
k}^{\left(  m\right)  }$ to $\mathcal{P}_{1\times k}^{\left(  m-2\right)  }$
and $H$ \emph{multiplies} \textup{(}elementwise\textup{) }$\mathcal{P}%
_{1\times k}^{\left(  m\right)  }$ by the \emph{number} $\left(  k/2\right)
+m$. An easy computation shows that%
\begin{equation}
\left[  E,X^{+}\right]  =2X^{+},\qquad\left[  E,X^{-}\right]  =-2X^{-}%
,\qquad\left[  X^{-},X^{+}\right]  =E. \label{eqFD.19}%
\end{equation}
Eq.\ \textup{(\ref{eqFD.19})} gives a faithful representation of the Lie
algebra $\mathrm{sl}_{2}\left(  \mathbb{R}\right)  $. Thus the dual action of
$H$ is given by this representation. The integrated form of this Lie algebra
representation is more subtle to describe: it is the \emph{metaplectic
representation of the two-sheeted covering group} $\widetilde{\mathrm{SL}%
_{2}\left(  \mathbb{R}\right)  }$ of $\mathrm{SL}_{2}\left(  \mathbb{R}%
\right)  $ \textup{(}or $\mathrm{Sp}_{2}\left(  \mathbb{R}\right)
$\textup{),} and this group is not a matrix group. Its concrete description
can be obtained by applying the Bargmann--Segal transform which sends the
Schr\"{o}dinger representation of this group to its Fock representation
$\mathcal{F}_{1\times k}$. However, for our purpose, its infinitesimal action
\textup{(\ref{eqFD.19})} together with the action of its maximal compact group
$G^{\prime}=\mathrm{U}\left(  1\right)  $, which is particularly simple, will
suffice. Indeed, it is easy to show that we have the following decomposition
of $\mathcal{P}_{1\times k}^{\left(  m\right)  }$:%
\begin{equation}
\mathcal{P}_{1\times k}^{\left(  m\right)  }=\mkern27mu\sum\limits_{\makebox
[0pt]{\hss$\scriptstyle i=0,\dots,\left[ m/2\right] $\hss}}{\!{\oplus}%
\mkern9mu}p_{0}^{i}\mathcal{H}_{1\times k}^{\left(  m-2i\right)  },
\label{eqFD.20}%
\end{equation}
where $\left[  m/2\right]  $ denotes the integral part of $m/2$, and
$\mathcal{H}_{1\times k}^{\left(  m-2i\right)  }$ denotes the subspace of all
harmonic homogeneous polynomials of degree $\left(  m-2i\right)  $, i.e., all
functions $p\in\mathcal{P}_{1\times k}^{\left(  m-2i\right)  }$ such that
$\bigtriangleup p=0$. For an integer $r\geq0$ then it can be easily shown that
the restriction $R_{H}^{\left(  r\right)  }$ of $R_{H}$ to $\mathcal{H}%
_{1\times k}^{\left(  r\right)  }$ is an irreducible representation of $H$
with signature $(\underset{\left[  k/2\right]  }{\underbrace{r,0,\dots,0}})$
and highest weight vector
\begin{equation}
f^{\left(  r\right)  _{H}}\left(  Z\right)  =%
%TCIMACRO{\TeXButton{cases}{\begin{cases}
%\left( Z_{1}+iZ_{s+1}\right) ^{r}, & \text{if }k=2s,  \\
%\left( Z_{1}+iZ_{s+2}\right) ^{r}, & \text{if }k=2s+1,\qquad i=\sqrt{-1}.
%\end{cases}
%}}%
%BeginExpansion
\begin{cases}
\left( Z_{1}+iZ_{s+1}\right) ^{r}, & \text{if }k=2s,  \\
\left( Z_{1}+iZ_{s+2}\right) ^{r}, & \text{if }k=2s+1,\qquad i=\sqrt{-1}.
\end{cases}
%EndExpansion
\label{eqFD.21}%
\end{equation}
For each integer $j\geq0$, the restriction of $R_{H}$ to the subspace
$p_{0}^{j}\mathcal{H}^{\left(  r\right)  }$ is equivalent to $R_{H}^{\left(
r\right)  }$ since $p_{0}^{j}$ is $H$-invariant. Set
\begin{equation}
\mathcal{I}_{1\times k}^{\left(  r\right)  }=\sum\limits_{\makebox
[0pt]{\hss$\scriptstyle j=0$\hss}}^{\infty}{\!{\oplus}\,}p_{0}^{j}%
\mathcal{H}_{1\times k}^{\left(  r\right)  }; \label{eqFD.22}%
\end{equation}
then $\mathcal{I}_{1\times k}^{\left(  r\right)  }$ is the $(\underset
{k}{\underbrace{r,0,\dots,0}})$-isotypic component of $R_{H}^{\left(
r\right)  }$. From \textup{(\ref{eqFD.20})} and \textup{(\ref{eqFD.22})} we
see that%
\begin{equation}
\mathcal{F}_{1\times k}=\sum\limits_{\makebox[0pt]{\hss$\scriptstyle r=0$\hss
}}^{\infty}{\!{\oplus}\,}\mathcal{I}_{1\times k}^{\left(  r\right)  }.
\label{eqFD.23}%
\end{equation}
Obviously, $R_{H^{\prime}}^{\prime}\left(  u\right)  =R_{G}^{\prime}\left(
u\right)  $, $u\in G^{\prime}$, leaves each one-dimensional subspace
$cp_{0}^{j}h$, $c\in\mathbb{C}$, invariant, since $R_{G}^{\prime}\left(
u\right)  \left(  p_{0}^{j}h\right)  =u^{r+2j}\left(  p_{0}^{j}h\right)  $
\textup{(}alternatively, $E\left(  p_{0}^{j}h\right)  =\left(  \left(
k/2\right)  +r+2j\right)  \left(  p_{0}^{j}h\right)  $\textup{),} for all
$h\in\mathcal{H}_{1\times k}^{\left(  r\right)  }$. Clearly, $X^{+}\left(
p_{0}^{j}h\right)  =\frac{1}{2}p_{0}^{\left(  j+1\right)  }h$, $h\in
\mathcal{H}_{1\times k}^{\left(  r\right)  }$. Finally from the equation%
\begin{equation}
X^{-}\left(  p_{0}f\right)  =\left(  k+2s\right)  f+\frac{1}{2}p_{0}%
\bigtriangleup f \label{eqFD.24}%
\end{equation}
if $f$ is a polynomial function of degree $s$, we deduce by induction on the
integer $j\geq1$ that%
\begin{equation}
X^{-}\left(  p_{0}^{j}h\right)  =j\left(  k+2\left(  r+j-1\right)  \right)
p_{0}^{j-1}h,\qquad h\in\mathcal{H}_{1\times k}^{\left(  r\right)  }.
\label{eqFD.25}%
\end{equation}
For each fixed $h\in\mathcal{H}_{1\times k}^{\left(  r\right)  }$ let $Jh$
denote the subspace of $\mathcal{I}_{1\times k}^{\left(  r\right)  }$ spanned
by the set $\left\{  p_{0}^{j}h\mid j=0,1,2,\dots\right\}  $. Then it follows
from the previous discussion that the subrepresentation of the Lie algebra
$\mathrm{sl}_{2}\left(  \mathbb{R}\right)  $ on $Jh$ is irreducible, and thus
the metaplectic subrepresentation of $\widetilde{\mathrm{SL}_{2}\left(
\mathbb{R}\right)  }$ on $Jh$ is irreducible as well. As a $\mathrm{U}\left(
1\right)  $-module $Jh$ is reducible, and for this special case each
one-dimensional subspace $cp_{0}^{j}h$, $c\in\mathbb{C}$, is an irreducible
submodule, and the lowest one is $ch$ which has weight $r$ \textup{(}or
$\left(  k/2\right)  +r$\textup{)} since%
\begin{equation}
R_{G}^{\prime}\left(  u\right)  h=u^{r}h,\qquad u\in\mathrm{U}\left(
1\right)  ,\text{\qquad or\qquad}Eh=\left(  \frac{k}{2}+r\right)  h.
\label{eqFD.26}%
\end{equation}
In general, if a \emph{holomorphic discrete series} of a noncompact semisimple
Lie group such as $\widetilde{\mathrm{SL}_{2}\left(  \mathbb{R}\right)  }$
considered as a $K$-module, where $K$ is its maximal compact subgroup,
decomposes into a discrete sum of irreducible submodules, each one of them can
be characterized by a signature \textup{(}highest weight, for
example\textup{)} and the one with the lowest highest weight \textup{(}under
the lexicographic ordering\textup{)} is unique. This \emph{lowest}
$K$\emph{-type highest weight} which corresponds to the \emph{Harish Chandra's
or Blattner's parameter,} can be used to label the given holomorphic discrete
series. We shall call this label its \emph{signature.} In our example, the
holomorphic discrete series $Jh$ of $\widetilde{\mathrm{SL}_{2}\left(
\mathbb{R}\right)  }$ has signature $r$. If $\dim\left(  \mathcal{H}_{1\times
r}^{\left(  r\right)  }\right)  =d$ \textup{(}actually, $d=\binom{k+r-1}%
{r}-\binom{k+r-3}{r-2}$\textup{)} $\mathcal{I}_{1\times k}^{\left(  r\right)
}$ is the $r$-isotypic component \textup{(}of the metaplectic representation
of $\widetilde{\mathrm{SL}_{2}\left(  \mathbb{R}\right)  }$\textup{)} which
contains $d$ isomorphic copies of signature $r$.

Now let us verify Theorem \textup{\ref{ThmFD.2}} for this simple example. From
Eq.\ \textup{(\ref{eqFD.20}) we have}%
\begin{multline}
\dim\left[  \operatorname*{Hom}\nolimits_{\mathrm{SO}\left(  k\right)
}\left(  (\underset{\left[  k/2\right]  }{\underbrace{r,0,\dots,0}%
})_{\mathrm{SO}\left(  k\right)  }:(\underset{k}{\underbrace{m,0,\dots,0}%
})_{\mathrm{U}\left(  k\right)  }\bigg|_{\mathrm{SO}\left(  k\right)
}\right)  \right] \label{eqFD.27}\\
=%
%TCIMACRO{\TeXButton{cases}{\begin{cases}
%1, & \text{if }r=m-2i\text{ for }i=0,\dots,\left[ m/2\right] ,  \\
%0, & \text{otherwise,}
%\end{cases}
%}}%
%BeginExpansion
\begin{cases}
1, & \text{if }r=m-2i\text{ for }i=0,\dots,\left[ m/2\right] ,  \\
0, & \text{otherwise,}
\end{cases}
%EndExpansion
\end{multline}
and from Eq.\ \textup{(\ref{eqFD.22})} and Eq.\ \textup{(\ref{eqFD.26})} we
have%
\begin{equation}
\dim\left[  \operatorname*{Hom}\nolimits_{\mathrm{U}\left(  1\right)  }\left(
m_{\mathrm{U}\left(  1\right)  }:r_{\widetilde{\mathrm{SL}_{2}\left(
\mathbb{R}\right)  }}\bigg|_{\mathrm{U}\left(  1\right)  }\right)  \right]  =%
%TCIMACRO{\TeXButton{cases}{\begin{cases}
%1, & \text{if }2j+r=m,  \\
%0, & \text{otherwise,}
%\end{cases}
%}}%
%BeginExpansion
\begin{cases}
1, & \text{if }2j+r=m,  \\
0, & \text{otherwise,}
\end{cases}
%EndExpansion
\label{eqFD.28}%
\end{equation}
which are obviously identical.

For arbitrary $n$ such that $n\leq k$ Eq.\ \textup{(\ref{eqFD.7})} remains
valid with $\left(  r\right)  =\left(  r_{11},\dots,r_{nk}\right)  $ and
$Z^{\left(  r\right)  }=Z_{11}^{r_{11}}\cdots Z_{nk}^{r_{nk}}$.
Eq.\ \textup{(\ref{eqFD.8}), (\ref{eqFD.12}), (\ref{eqFD.13}), (\ref{eqFD.14})
remain valid. Eq.\ (\ref{eqFD.10}) is replaced by}%
\begin{equation}
R_{ij}=\sum_{\alpha=1}^{n}Z_{\alpha i}\frac{\partial\;}{\partial Z_{\alpha j}%
},\qquad1\leq i,j\leq k. \tag*{$($\textup{\ref{eqFD.10}}$)^{\prime}$}%
\end{equation}
Eq.\ \textup{(\ref{eqFD.11})} is replaced by%
\begin{equation}
L_{\alpha\beta}=\sum_{i=1}^{k}Z_{\alpha i}\frac{\partial\;}{\partial Z_{\beta
i}},\qquad1\leq\alpha,\beta\leq n. \tag*{$($\textup{\ref{eqFD.11}}$)^{\prime
}$}%
\end{equation}
Let $B_{n}^{\prime}$ denote the lower triangular Borel subgroup of
$G_{\mathbb{C}}^{\prime}=\mathrm{GL}_{n}\left(  \mathbb{C}\right)  $, let
$\left(  \lambda\right)  $ be an $n$-tuple of integers such that $\lambda
_{1}\geq\lambda_{2}\geq\dots\geq\lambda_{n}\geq0$, let $\lambda\colon
B_{n}^{\prime}\rightarrow\mathbb{C}^{\ast}$ be the holomorphic character
defined on $B_{n}^{\prime}$ by%
\[
\lambda\left(  b^{\prime}\right)  =\left(  b_{11}^{\prime}\right)
^{\lambda_{1}}\cdots\left(  b_{nn}^{\prime}\right)  ^{\lambda_{n}}\text{\qquad
if }b^{\prime}=%
\begin{bmatrix}
b_{11}^{\prime} &  & \raisebox{-9pt}[0pt][0pt]{\llap{\huge{$0$}\kern3pt}}\\
& \ddots & \\
\raisebox{3pt}[0pt][0pt]{\rlap{\kern3pt\huge{$\ast$}}} &  & b_{nn}^{\prime}%
\end{bmatrix}
\text{ belongs to }B_{n}^{\prime}.
\]
Let $\mathcal{P}_{n\times k}^{\left(  \lambda\right)  }$ denote the subspace
of all polynomial functions on $\mathbb{C}^{n\times k}$ which also satisfy the
covariant condition%
\begin{equation}
f\left(  b^{\prime}Z\right)  =\lambda\left(  b^{\prime}\right)  f\left(
Z\right)  ,\qquad\left(  b^{\prime},Z\right)  \in B_{n}^{\prime}%
\times\mathbb{C}^{n\times k}. \label{eqFD.29}%
\end{equation}
Let $R_{\lambda}$ denote the representation of $G$ obtained by right
translation on $\mathcal{P}_{n\times k}^{\left(  \lambda\right)  }$. Then by
the Borel--Weil theorem \textup{(}see, e.g., \cite[Theorem 1.5]{Ton76a}%
\textup{)} $R_{\lambda}$ is irreducible with highest weight $\left(
\lambda\right)  $ and highest weight vector%
\begin{equation}
cf_{\lambda}\left(  Z\right)  =c\Delta_{1}^{\lambda_{1}-\lambda_{2}}\left(
Z\right)  \Delta_{2}^{\lambda_{2}-\lambda_{3}}\left(  Z\right)  \cdots
\Delta_{{}}^{\lambda_{n}}\left(  Z\right)  ,\qquad c\in\mathbb{C}^{\ast},
\label{eqFD.30}%
\end{equation}
where in Eq.\ \textup{(\ref{eqFD.30})} $\Delta_{i}\left(  Z\right)  $ denotes
the $i^{\text{th}}$ principal minor of $Z$.

Similarly let $B_{k}^{t}$ denote the upper triangular Borel subgroup of
$G_{\mathbb{C}}=\mathrm{GL}_{k}\left(  \mathbb{C}\right)  $ and let
$\lambda^{\prime}\colon B_{k}^{t}\rightarrow\mathbb{C}^{\ast}$ be the
holomorphic character defined on $B_{k}^{t}$ by
\[
\lambda^{\prime}\left(  b\right)  =b_{11}^{\lambda_{1}}\cdots b_{nn}%
^{\lambda_{n}}\text{\qquad if }b=%
\begin{bmatrix}
b_{11} &  &  &  & \raisebox{-12pt}[0pt][0pt]{\llap{\Huge{$\ast$}\kern6pt}}\\
& \ddots &  &  & \\
&  &  b_{nn} &  & \\
&  &  & \ddots & \\
\raisebox{6pt}[0pt][0pt]{\rlap{\kern6pt\Huge{$0$}}} &  &  &  & b_{kk}%
\end{bmatrix}
\text{ belongs to }B_{k}^{t}.
\]

Let $\mathcal{P}_{n\times k}^{\left(  \lambda^{\prime}\right)  }$ denote the
subspace of all polynomial functions on $\mathbb{C}^{n\times k}$ which also
satisfy the covariant condition%
\begin{equation}
f\left(  Zb\right)  =\lambda^{\prime}\left(  b\right)  f\left(  Z\right)
,\qquad\left(  b,Z\right)  =B_{k}^{t}\times\mathbb{C}^{n\times k}\text{.}
\label{eqFD.31}%
\end{equation}
Let $R_{\lambda^{\prime}}^{\prime}$ denote the representation of $G^{\prime}$
on $\mathcal{P}_{n\times k}^{\left(  \lambda^{\prime}\right)  }$ defined by%
\begin{equation}
\left[  R_{\lambda^{\prime}}^{\prime}\left(  g^{\prime}\right)  f\right]
\left(  Z\right)  =f\left(  \left(  g^{\prime}\right)  ^{t}Z\right)  ,\qquad
g^{\prime}\in G^{\prime}. \label{eqFD.32}%
\end{equation}
Then $R_{\lambda^{\prime}}^{\prime}$ is irreducible with highest weight
$\left(  \lambda^{\prime}\right)  $ and with the same highest weight vector
given by Eq.\ \textup{(\ref{eqFD.30}).} By Weyl's unitarian trick the
restriction of $R_{\lambda}$ \textup{(}resp.\ $R_{\lambda^{\prime}}^{\prime}%
$\textup{)} to $G=\mathrm{U}\left(  k\right)  $ \textup{(}resp.\ $G^{\prime
}=\mathrm{U}\left(  n\right)  $\textup{)} remains irreducible with the same 
signature.

Let $\mathcal{I}_{n\times k}^{\left(  \lambda\right)  }$ denote the
$G_{\mathbb{C}}^{\prime}\times G_{\mathbb{C}}^{{}}$ \textup{(}or $G^{\prime
}\times G$\textup{)}-cyclic module in $\mathcal{F}_{n\times k}$ generated by
the highest vector $f_{\lambda}$ given by Eq.\ \textup{(\ref{eqFD.29});} then
by Theorem \textup{3,} p.\ \textup{150,} of \cite{Zel73}, $\mathcal{I}%
_{n\times k}^{\left(  \lambda\right)  }$ is irreducible with highest weight
$\left(  \lambda^{\prime},\lambda\right)  $. For the sake of simplicity we say
that the $G_{\mathbb{C}}^{\prime}\times G_{\mathbb{C}}^{{}}$-module
$\mathcal{I}_{n\times k}^{\left(  \lambda\right)  }$ has \emph{signature}
$\left(  \lambda\right)  $. To prove that $\mathcal{I}_{n\times k}^{\left(
\lambda\right)  }\approx\mathcal{P}_{n\times k}^{\left(  \lambda^{\prime
}\right)  }\mathop{\hat{\otimes}}\mathcal{P}_{n\times k}^{\left(
\lambda\right)  }$ we define a map $\Phi\colon\mathcal{P}_{n\times k}^{\left(
\lambda^{\prime}\right)  }\mathop{\hat{\otimes}}\mathcal{P}_{n\times
k}^{\left(  \lambda\right)  }\rightarrow\mathcal{I}_{n\times k}^{\left(
\lambda\right)  }$ as follows:

Let $f^{\prime}\otimes f\in\mathcal{P}_{n\times k}^{\left(  \lambda^{\prime
}\right)  }\mathop{\hat{\otimes}}\mathcal{P}_{n\times k}^{\left(
\lambda\right)  }$. Then $f^{\prime}$ and $f$ can be represented in the
following form:%
\begin{equation}
f^{\prime}=\sum_{i\in I^{\prime}}c_{i}^{\prime}R_{\lambda^{\prime}}^{\prime
}\left(  g_{i}^{\prime}\right)  f_{\lambda}^{{}},\qquad f=\sum_{j\in I}%
c_{j}R_{\lambda}\left(  g_{j}\right)  f_{\lambda}, \label{eqFD.33}%
\end{equation}
where in Eq.\ \textup{(\ref{eqFD.33})} $c_{i}^{\prime},c_{j}^{{}}\in
\mathbb{C}$, $g_{i}^{\prime}\in G_{\mathbb{C}}^{\prime}$, $g_{j}\in
G_{\mathbb{C}}$, and $I^{\prime}$ and $I$ are two finite index sets. Set
$\Phi\left(  f^{\prime}\otimes f\right)  =\sum_{i\in I^{\prime},\,j\in I}%
c_{i}^{\prime}c_{j}^{{}}T\left(  g_{i}^{\prime},g_{j}^{{}}\right)  f_{\lambda
}^{{}}$, where $\left[  T\left(  g_{i}^{\prime},g_{j}^{{}}\right)  f_{\lambda
}\right]  \left(  Z\right)  =f\left(  \left(  g_{i}^{\prime}\right)
^{t}Zg_{j}^{{}}\right)  $. Since
\[
R_{\lambda^{\prime}}^{\prime}\left(  g^{\prime}\right)  f^{\prime}=\sum_{i\in
I}c_{i}^{\prime}R_{\lambda^{\prime}}^{\prime}\left(  g_{{}}^{\prime}%
g_{i}^{\prime}\right)  f_{\lambda}^{{}}%
\]
and
\[
R_{\lambda}\left(  g\right)  f=\sum_{j\in I}c_{j}R_{\lambda}\left(
gg_{j}\right)  f_{\lambda}%
\]
it follows that%
\begin{align*}
\Phi\left[  \left(  R_{\lambda^{\prime}}^{\prime}\left(  g^{\prime}\right)
\otimes R_{\lambda}\left(  g\right)  \right)  \left(  f^{\prime}\otimes
f\right)  \right]   &  =\sum_{i\in I,\,j\in T}c_{i}^{\prime}c_{j}^{{}}T\left(
g_{{}}^{\prime}g_{i}^{\prime},gg_{j}^{{}}\right)  f_{\lambda}^{{}}\\
&  =T\left(  g^{\prime},g\right)  \Phi\left(  f^{\prime}\otimes f\right)
\end{align*}
for all $g^{\prime}\in G_{\mathbb{C}}^{\prime}$ and $g\in G_{\mathbb{C}}$.
This means that $\Phi$ is an intertwining operator and by Schur's lemma $\Phi$
is either $0$ or an isomorphism. Since%
\[
\Phi\left(  f_{\lambda}\otimes f_{\lambda}\right)  =f_{\lambda}%
\]
it follows that $\Phi$ is an isomorphism. Since $\mathcal{P}_{n\times k}$ is
dense in $\mathcal{F}_{n\times k}$ Theorem \textup{3} \textup{(}%
p.\ \textup{150)} of \cite{Zel73} \textup{(}see also \cite{KlTo89}\textup{)}
implies that we have the Hilbert sum $\mathcal{F}_{n\times k}= \smash
[b]{\sum\limits_{\left(  \lambda\right)  }{\!{\oplus}\,}\mathcal{I}_{n\times
k}^{\left(  \lambda\right)  }}$ for the pair $\left(  \mathrm{U}\left(
n\right)  ,\mathrm{U}\left(  k\right)  \right)  $.

Now suppose $k>2n$ and set $H=\mathrm{SO}\left(  k\right)  $, $H_{\mathbb{C}%
}=\mathrm{SO}_{k}\left(  \mathbb{C}\right)  $. Let $J_{n\times k}$ denote the
ring of all $H$ \textup{(}or $H_{\mathbb{C}}$\textup{)}-invariant polynomials
in $\mathcal{P}_{n\times k}$. Then $J_{n\times k}$ is generated by the
constants and the $n\left(  n+1\right)  /2$ algebraically independent
polynomials%
\begin{equation}
p_{\alpha\beta}\left(  Z\right)  =\sum_{i=1}^{k}Z_{\alpha i}Z_{\beta i}%
,\qquad1\leq\alpha\leq\beta\leq n. \label{eqFD.34}%
\end{equation}
It follows that the ring of all $H$ \textup{(}or $H_{\mathbb{C}}$%
\textup{)}-invariant differential operators with constant coefficients is
generated by the constants and the Laplacians%
\begin{equation}
\bigtriangleup_{\alpha\beta}=p_{\alpha\beta}\left(  D\right)  =\sum_{i=1}%
^{k}\frac{\partial^{2}\;}{\partial Z_{\alpha i}\partial Z_{\beta i}}%
,\qquad1\leq\alpha\leq\beta\leq n. \label{eqFD.35}%
\end{equation}
The infinitesimal action of $R_{G_{\mathbb{C}}}^{\prime}$ is generated by
\begin{equation}
L_{\alpha\beta}=\sum_{i=1}^{k}Z_{\alpha i}\frac{\partial\;}{\partial Z_{\beta
i}},\qquad1\leq\alpha,\beta\leq n. \label{eqFD.36}%
\end{equation}
Set $P_{\alpha\beta}=-p_{\alpha\beta}$, $E_{\alpha\beta}=L_{\alpha\beta}%
+\frac{1}{2}k\delta_{\alpha\beta}$, and $D_{\alpha\beta}=\bigtriangleup
_{\alpha\beta}$; then it follows from \cite{KLT95} \textup{(}see
Eq.\ \textup{(3.3))} that $\left\{  E_{\alpha\beta},P_{\alpha\beta}%
,D_{\alpha\beta}\right\}  $ defines a faithful representation of
$\mathrm{sp}_{2n}\left(  \mathbb{R}\right)  $ on $\mathcal{F}_{n\times k}$. By
construction this representation is dual to the infinitesimal action of
$R_{H}$. The global action $R_{H^{\prime}}^{\prime}$ is a unitary metaplectic
representation of $\widetilde{\mathrm{Sp}_{2n}\left(  \mathbb{R}\right)  }$,
the two-sheeted covering of $\mathrm{Sp}_{2n}\left(  \mathbb{R}\right)  $
\textup{(}see \cite{KLT95} for details\textup{).} As in the case of the pair
$\left(  \mathrm{U}\left(  n\right)  ,\mathrm{U}\left(  k\right)  \right)  $
the common highest weight vector \textup{(}for $R_{H^{\prime}}^{\prime}$ the
lowest $K^{\prime}$-type highest weight vector\textup{)} of signature $\left(
\mu\right)  =\left(  \mu^{}_{1},\dots,\mu^{}_{n}\right)  $ with $\mu^{}%
_{1}\geq\dots\geq\mu^{}_{n}\geq0$ and $\mu^{}_{i}\in\mathbb{N}$, $1\leq i\leq
n$, of the pair $\left(  \widetilde{\mathrm{Sp}_{2n}\left(  \mathbb{R}\right)
,}\mathrm{SO}\left(  k\right)  \right)  $ is%
\begin{equation}
f_{\mu}\left(  Z\right)  =\Delta_{1}^{\mu^{}_{1}-\mu^{}_{2}}\left(  Zq\right)
\Delta_{2}^{\mu^{}_{2}-\mu^{}_{3}}\left(  Zq\right)  \cdots\Delta_{n}^{\mu
^{}_{n}}\left(  Zq\right)  , \label{eqFD.37}%
\end{equation}
where the $k\times k$ matrix $q$ is given by
\[
\frac{1}{\sqrt{2}}\left[
\begin{tabular}
[c]{c|c}%
$\openone_{\nu_{\mathstrut}}$ & $\openone_{\nu_{\mathstrut}}$\\\hline
$i\openone_{\nu_{\mathstrut}}$ & $-i\openone_{\nu_{\mathstrut}}$%
\end{tabular}
\right]  \text{\qquad if }k=2\nu,
\]
and
\[
\frac{1}{\sqrt{2}}%
\begin{bmatrix}
\openone_{\nu} & 0 & \openone_{\nu}\\
0 & \sqrt{2} & 0\\
i\openone_{\nu} & 0 & -i\openone_{\nu}%
\end{bmatrix}
\text{\qquad if }k=2\nu+1,
\]
and where $\openone_{\nu}$ is the unit matrix of order $\nu$.

An element $p$ of $\mathcal{P}_{n\times k}$ is called $H$\emph{-harmonic} if
$\bigtriangleup_{\alpha\beta}p=0$ for all $\alpha,\beta=1,\dots,n$. Let
$\mathcal{H}_{n\times k}$ denote the subspace of all $H$-harmonic polynomial
functions of $\mathcal{P}_{n\times k}$ and let $\mathcal{H}_{n\times k}\left(
\mu\right)  $ denote the subspace of all elements $h$ of $\mathcal{H}_{n\times
k}$ which also satisfy the covariant condition%
\begin{equation}
h\left(  b^{\prime}Z\right)  =\left(  b_{11}^{\prime}\right)  ^{\mu^{}_{1}%
}\cdots\left(  b_{nn}^{\prime}\right)  ^{\mu^{}_{n}}h\left(  Z\right)
,\qquad\forall\,b^{\prime}\in B_{n}^{\prime}. \label{eqFD.38}%
\end{equation}
Then according to Theorem \textup{3.1} of \cite{Ton76a}, the representation
$R_{H}$ of $H$ which is obtained by right translations on $\mathcal{H}%
_{n\times k}\left(  \mu\right)  $ is irreducible with signature $\left(
\mu\right)  $.

The infinitesimal action of $R_{H}$ is given by%
\begin{equation}
R_{ij}^{H}=\sum_{\alpha=1,\dots,n}\left(  Z_{\alpha i}\frac{\partial
\;}{\partial Z_{\alpha j}}-Z_{\alpha j}\frac{\partial\;}{\partial Z_{\alpha
i}}\right)  ,\qquad1\leq i<j\leq k. \label{eqFD.39}%
\end{equation}
{}From \cite{KLT95} the dual infinitesimal action of $R_{H}$ is given by the
system $\left\{  E_{\alpha\beta},P_{\alpha\beta},D_{\alpha\beta}\right\}  $
which satisfies the commutation relations%
\begin{equation}%
%TCIMACRO{\TeXButton{cases}{\begin{cases}
%\left[  E_{\alpha\beta},E_{\mu\nu}\right]     =\delta_{\beta\mu}E_{\alpha\nu
%}-\delta_{\alpha\nu}E_{\mu\beta}\\
%\left[  E_{\alpha\beta},P_{\mu\nu}\right]     =\delta_{\beta\mu}P_{\alpha\nu
%}+\delta_{\beta\nu}P_{\alpha\mu}\\
%\left[  E_{\alpha\beta},D_{\mu\nu}\right]     =-\delta_{\alpha\mu}D_{\beta
%\nu}-\delta_{\alpha\nu}D_{\beta\mu}\\
%\left[  P_{\alpha\beta},D_{\mu\nu}\right]     =\delta_{\alpha\mu}E_{\nu\beta
%}+\delta_{\alpha\nu}E_{\mu\beta}+\delta_{\beta\mu}E_{\nu\alpha}+\delta
%_{\beta\nu}E_{\mu\alpha}\\
%\left[  P_{\alpha\beta},P_{\mu\nu}\right]     =\left[  D_{\alpha\beta}%
%,D_{\mu\nu}\right]  =0\\
%P_{\alpha\beta},\;P_{\beta\alpha},\;D_{\alpha\beta}   =D_{\beta\alpha}\\
%P_{\alpha\beta}^{\dag}   =D_{\alpha\beta}^{{}},\qquad D_{\alpha\beta}^{\dag
%}=P_{\alpha\beta}^{{}},\qquad E_{\alpha\beta}^{\dag}=E_{\beta\alpha}^{{}},  \\
%\text{\qquad for all }\alpha,\beta,\mu,\nu=1,\dots,n.
%\end{cases}}}%
%BeginExpansion
\begin{cases}
\left[  E_{\alpha\beta},E_{\mu\nu}\right]     =\delta_{\beta\mu}E_{\alpha\nu
}-\delta_{\alpha\nu}E_{\mu\beta}\\
\left[  E_{\alpha\beta},P_{\mu\nu}\right]     =\delta_{\beta\mu}P_{\alpha\nu
}+\delta_{\beta\nu}P_{\alpha\mu}\\
\left[  E_{\alpha\beta},D_{\mu\nu}\right]     =-\delta_{\alpha\mu}D_{\beta
\nu}-\delta_{\alpha\nu}D_{\beta\mu}\\
\left[  P_{\alpha\beta},D_{\mu\nu}\right]     =\delta_{\alpha\mu}E_{\nu\beta
}+\delta_{\alpha\nu}E_{\mu\beta}+\delta_{\beta\mu}E_{\nu\alpha}+\delta
_{\beta\nu}E_{\mu\alpha}\\
\left[  P_{\alpha\beta},P_{\mu\nu}\right]     =\left[  D_{\alpha\beta}%
,D_{\mu\nu}\right]  =0\\
P_{\alpha\beta},\;P_{\beta\alpha},\;D_{\alpha\beta}   =D_{\beta\alpha}\\
P_{\alpha\beta}^{\dag}   =D_{\alpha\beta}^{{}},\qquad D_{\alpha\beta}^{\dag
}=P_{\alpha\beta}^{{}},\qquad E_{\alpha\beta}^{\dag}=E_{\beta\alpha}^{{}},  \\
\text{\qquad for all }\alpha,\beta,\mu,\nu=1,\dots,n.
\end{cases}%
%EndExpansion
\label{eqFD.40}%
\end{equation}
By Corollary \textup{3.11} of \cite{Ton76a} the $\mu$-isotypic component in
$\mathcal{H}_{n\times k}$ consists of $d_{\mu}$ copies isomorphic to
$\mathcal{H}_{n\times k}\left(  \mu\right)  $, where $d_{\mu}$ is the degree
of an irreducible representation of $G^{\prime}=\mathrm{U}\left(  n\right)  $
of signature $\left(  \mu^{}_{1},\dots,\mu^{}_{n}\right)  $. Since from
Eq.\ \textup{(\ref{eqFD.40})} and the fact that $f_{\mu}$ is $H$-harmonic%
\begin{align*}
D_{\mu\nu}E_{\alpha\beta}f_{\mu}  &  =\left[  D_{\mu\nu},E_{\alpha\beta
}\right]  f_{\mu}+E_{\alpha\beta}D_{\mu\nu}f_{\mu}\\
&  =\delta_{\alpha\mu}D_{\beta\nu}f_{\mu}+\delta_{\alpha\nu}D_{\beta\mu}%
f_{\mu}\\
&  =0,
\end{align*}
it follows that $E_{\alpha\beta}f_{\mu}$ is $H$-harmonic for every
$\alpha,\beta=1,\dots,n$. Since $\left[  E_{\alpha\beta}^{{}},R_{ij}%
^{H}\right]  =0$ for all $\alpha,\beta=1,\dots,n$ and $i,j=1,\dots,k$ it
follows that $E_{\alpha\beta}\colon\mathcal{H}_{n\times k}\left(  \mu\right)
\rightarrow\mathcal{H}_{n\times k}$ are intertwining operators, and thus are
either $0$ or isomorphisms. It follows that the $\mathfrak{g}^{\prime}$-module
generated by the cyclic vector $f_{\mu}$ is irreducible with signature
$\left(  \mu_{1},\dots,\mu^{}_{n}\right)  $. In fact, from
Eq.\ \textup{(3.14)} of \cite{Ton76a} this space is a $G^{\prime}$-module. Let
$G^{\prime}f_{\mu}$ denote this $G^{\prime}$-module; then by construction
$G^{\prime}f_{\mu}\subset\mathcal{H}_{n\times k}$.

If $h\in G^{\prime}f_{\mu}$ then from Eq.\ \textup{(\ref{eqFD.40})} we have%
\begin{align*}
D_{\mu\nu}P_{\alpha\beta}h  &  =\left[  D_{\mu\nu},P_{\alpha\beta}\right]
h+P_{\alpha\beta}D_{\mu\nu}h\\
&  =-\left(  \delta_{\alpha\mu}E_{\nu\beta}+\delta_{\alpha\nu}E_{\mu\beta
}+\delta_{\beta\mu}E_{\nu\alpha}+\delta_{\beta\nu}E_{\mu\nu}\right)  h,
\end{align*}
and therefore $D_{\mu\nu}P_{\alpha\beta}h$ belongs to $G^{\prime}f_{\mu}$. It
follows that $J_{n\times k}G^{\prime}f_{\mu}$ is an irreducible $\mathrm{sp}%
_{2n}\left(  \mathbb{R}\right)  $-module with signature $\left(  \mu\right)
$. Let $\mathcal{H}_{n\times k}^{\prime}\left(  \mu\right)  $ denote this
module and let $\mathcal{I}_{n\times k}^{\left(  \mu\right)  }$ be the
$H^{\prime}\times H$-cyclic module generated by $f_{\mu}$; then a proof
similar to the case $\mathcal{I}_{n\times k}^{\left(  \lambda\right)  }$ shows
that $\mathcal{H}_{n\times k}^{\prime}\left(  \mu\right)  \mathop{\hat
{\otimes}} \mathcal{H}_{n\times k}^{{}}\left(  \mu\right)  $ is isomorphic to
$\mathcal{I}_{n\times k}^{\left(  \mu\right)  }$. By the ``separation of
variables theorem'' \textup{2.5} of \cite{Ton76a} and from the fact that
$\mathcal{P}_{n\times k}$ is dense in $\mathcal{F}_{n\times k}$ it follows
that the orthogonal direct sum decomposition $\mathcal{F}_{n\times k}%
=\sum\limits_{\left(  \mu\right)  }{\!{\oplus}\,}\mathcal{I}_{n\times
k}^{\left(  \mu\right)  }$ holds. Therefore the reciprocity theorem
\textup{\ref{ThmFD.2}} holds for these pairs $\left(  G^{\prime},G\right)  $
and $\left(  H^{\prime},H\right)  $ as well. \medskip

%\item \label{ExaFD.3(2)} 
%
%
%
%
%
%
\noindent2) Let $k=2l$ and consider again the dual pair $\left(  G^{\prime
}=\mathrm{U}\left(  n\right)  ,G=\mathrm{U}\left(  k\right)  \right)  $. Let
$H=$\textrm{Sp}$\left(  k\right)  $; then $H_{\mathbb{C}}=$\textrm{Sp}%
$_{k}\left(  \mathbb{C}\right)  $. If $l\geq n\geq2$ then the theory of
symplectic harmonic polynomials in \cite{Ton77} implies that the dual
representation to the representation $R_{H}$ on $\mathcal{F}_{n\times k}$ is a
representation of the group $\mathrm{SO}^{\ast}\left(  2n\right)  =H^{\prime}$
whose infinitesimal action is given by Eq.\ \textup{(4.2)} of \cite{KLT95}.
Using Theorem \textup{2.1} of \cite{Ton77} and the ``separation of variables
theorem'' for this case we can show similarly that $\mathcal{F}_{n\times
k}=\sum\limits_{\left(  \mu\right)  }{\!{\oplus}\,}\mathcal{I}_{n\times
k}^{\left(  \mu\right)  }$ for this dual pair $\left(  \mathrm{SO}^{\ast
}\left(  2n\right)  ,\mathrm{Sp}\left(  k\right)  \right)  $. Thus the
reciprocity theorem \textup{\ref{ThmFD.2}} holds again for these pairs
$\left(  G^{\prime},G\right)  $ and $\left(  H^{\prime},H\right)  $. \medskip

%\item \label{ExaFD.3(3)} 
%
%
%
%
%
%
\noindent3) The case of the dual pairs $\left(  G^{\prime}=\mathrm{U}\left(
p\right)  \times\mathrm{U}\left(  q\right)  ,G=\mathrm{U}\left(  k\right)
\times\mathrm{U}\left(  k\right)  \right)  $ and \linebreak $\left(
H^{\prime}=\mathrm{U}\left(  p,q\right)  ,H=\mathrm{U}\left(  k\right)
\right)  $ can be treated in a similar fashion using the results of
\cite{Ton76b} and the infinitesimal action of $H^{\prime}$ on $\mathcal{F}%
_{n\times k}$ is given by Eq.\ \textup{(6.4)} of \cite{Ton95}. However, its
generalization to the case $H=\mathrm{U}\left(  \infty\right)  $ in Section
\textup{\ref{FID}} is quite delicate and requires a quite different embedding
that we shall describe in detail below.

Let $p$ and $q$ be positive integers such that $p+q=n$. Let $k$ be an integer
such that $k\geq2\max\left(  p,q\right)  $. Let $\left(  \lambda\right)  $ be
a $q$-tuple of integers such that $\lambda_{1}\geq\lambda_{2}\geq\dots
\geq\lambda_{q}\geq0$. Let $R_{\lambda}$ denote the representation of
$\mathrm{GL}_{k}\left(  \mathbb{C}\right)  $ \textup{(}or $\mathrm{U}\left(
k\right)  $\textup{)} defined on $\mathcal{P}_{q\times k}^{\left(
\lambda\right)  }$ given by Eq.\ \textup{(\ref{eqFD.29})} and
\textup{(\ref{eqFD.30})} with $n$ replaced by $q$. We define the
\emph{contragredient} \textup{(}or \emph{dual\/}\textup{)}
\emph{representation} of $R_{\lambda}$ as follows.

Let $s_{r}$ denote the $r\times r$ matrix with ones along the reverse diagonal
and zero elsewhere:%
\[%
\begin{pmatrix}
\raisebox{-7pt}[0pt][0pt]{\rlap{\kern3pt\huge{$0$}}} &  & 1\\
& \begin{picture}(0.5,1)(0.25,0.167) \multiput(1.25,1.25)(-0.25,-0.25){7}%
{\makebox(0,0){$\cdot$}} \end{picture} & \\
1 &  & \raisebox{1pt}[0pt][0pt]{\llap{\huge{$0$}\kern3pt}}%
\end{pmatrix}
.
\]
If $W\in\mathbb{C}^{q\times k}$ let $\tilde{W}=s_{q}Ws_{k}$. Thus $\tilde{W}$
is of the form%
\begin{equation}
\tilde{W}=%
\begin{bmatrix}
W_{q,k} & \cdots &  W_{q,1}\\
\vdots &  & \vdots\\
W_{1,k} & \cdots &  W_{1,1}%
\end{bmatrix}
. \label{eqFD.41}%
\end{equation}
Let $\mathcal{P}_{q\times k}^{\left(  \lambda^{\checkmark}\right)  }$ denote
the subspace of all polynomial functions in $\tilde{W}$ which also satisfy the
covariant condition%
\begin{equation}
f\left(  \tilde{b}^{\prime}\tilde{W}\right)  =\lambda\left(  b^{\prime
}\right)  f\left(  \tilde{W}\right)  \label{eqFD.42}%
\end{equation}
for all $b^{\prime}\in B_{q}^{\prime}$, where $B_{q}^{\prime}$ is the lower
triangular Borel subgroup of $\mathrm{GL}_{q}\left(  \mathbb{C}\right)  $, and
$\tilde{b}^{\prime}=s_{q}b^{\prime}s_{q}$.

Define the representation $R_{\lambda^{\checkmark}}$ of $\mathrm{GL}%
_{k}\left(  \mathbb{C}\right)  $ \textup{(}or $\mathrm{U}\left(  k\right)
$\textup{)} on $\mathcal{P}_{q\times k}^{\left(  \lambda^{\checkmark}\right)
}$ by%
\begin{equation}
\left[  R_{\lambda^{\checkmark}}\left(  g\right)  f\right]  \left(  \tilde
{W}\right)  =f\left(  \tilde{W}s_{k}gs_{k}\right)  ,\qquad g\in\mathrm{GL}%
_{k}\left(  \mathbb{C}\right)  . \label{eqFD.43}%
\end{equation}
Then $R_{\lambda^{\checkmark}}$ is irreducible with signature $\left(
\smash[b]{\underset
{k}{\underbrace{0,\dots,0,-\lambda_{q},-\lambda_{q-1},\dots,-\lambda_{1}}}%
}\right)  $ and \emph{lowest weight vector}%
\begin{equation}
cf_{\lambda^{\checkmark}}\left(  \tilde{W}\right)  =\Delta_{1}^{\lambda
_{1}-\lambda_{2}}\left(  \tilde{w}\right)  \Delta_{2}^{\lambda_{2}-\lambda
_{3}}\left(  \tilde{w}\right)  \cdots\Delta_{{}}^{\lambda_{q}}\left(
\tilde{w}\right)  ,\qquad c\in\mathbb{C}^{\ast}, \label{eqFD.44}%
\end{equation}
of weight $\left(  -\lambda_{1},-\lambda_{2},\dots,-\lambda_{q},0,\dots
,0\right)  $.

Let $\mathcal{P}_{q\times k}^{\left(  \lambda^{\checkmark}\right)  ^{\prime}}$
denote the subspace of all polynomial functions in $\tilde{W}$ which also
satisfy the covariant condition
\begin{equation}
f\left(  \tilde{W}\tilde{b}\right)  =\lambda\left(  b\right)  f\left(
\tilde{W}\right)  , \label{eqFD.45}%
\end{equation}
where $\tilde{b}=s_{k}bs_{k}$, $b\in B_{k}^{t}$ \textup{(}it follows that
$\tilde{b}$ is a lower triangular matrix of the form $\tilde{b}=\left(
\begin{smallmatrix}
b_{kk} &  & \raisebox{-6pt}[0pt][0pt]{\llap{\Large{$0$}\kern2pt}}\\
& \raisebox{0pt}[10pt]{$\ddots$} & \\
\raisebox{2pt}[0pt][0pt]{\rlap{\kern2pt\Large{$\ast$}}} &  & b_{11}%
\end{smallmatrix}
\right)  $\textup{).} Let $R_{\left(  \lambda^{\checkmark}\right)  ^{\prime}%
}^{\prime}$ denote the representation of $\mathrm{GL}_{q}\left(
\mathbb{C}\right)  $ \textup{(}or of $G^{\prime}=\mathrm{U}\left(  q\right)
$\textup{)} on $\mathcal{P}_{q\times k}^{\left(  \lambda^{\checkmark}\right)
^{\prime}}$ defined by%
\begin{equation}
\left[  R_{\left(  \lambda^{\checkmark}\right)  ^{\prime}}^{\prime}\left(
g^{\prime}\right)  f\right]  \left(  \tilde{W}\right)  =f\left(  s_{q}\left(
g^{\prime}\right)  ^{-1}s_{q}\tilde{W}\right)  ,\qquad g^{\prime}%
\in\mathrm{GL}_{q}\left(  \mathbb{C}\right)  . \label{eqFD.46}%
\end{equation}
Then $R_{\left(  \lambda^{\checkmark}\right)  ^{\prime}}^{\prime}$ is
irreducible with highest weight $\left(  \lambda^{\checkmark}\right)
^{\prime}$ and with \emph{lowest weight vector} given by $cf_{\lambda
^{\checkmark}}$, $c\in\mathbb{C}^{\ast}$, of weight $\left(  -\lambda
_{1},-\lambda_{2},\dots,-\lambda_{q}\right)  $.

As in the case $\left(  \lambda^{\prime}\right)  \otimes\left(  \lambda
\right)  $ it can be shown that $\mathcal{P}_{q\times k}^{\left(
\lambda^{\checkmark}\right)  ^{\prime}}\mathop{\hat{\otimes}}\mathcal{P}%
_{q\times k}^{\left(  \lambda^{\checkmark}\right)  }$ is isomorphic to
$\mathcal{I}_{q\times k}^{\left(  \lambda^{\checkmark}\right)  }$ and we have
the Hilbert sum decomposition $\mathcal{F}_{q\times k}=\smash[b]{
\sum
\limits_{
\left(  \lambda^{\checkmark}\right)
}{\!{\oplus}\,}
\mathcal{I}_{q\times k}^{
\left(  \lambda^{\checkmark}\right)
}
}$ for the pair $\left(  \mathrm{U}\left(  q\right)  ,\mathrm{U}\left(
k\right)  \right)  $.

Now let $G=\mathrm{U}\left(  k\right)  \times\mathrm{U}\left(  k\right)  $ act
on $\mathcal{F}_{n\times k}$ via the outer tensor product%
\begin{equation}
\left[  R_{\mathrm{U}\left(  k\right)  }\mathop{\hat{\otimes}}R_{\mathrm{U}%
\left(  k\right)  ^{\checkmark}}\right]  \left(  g_{1},g_{2}\right)  f\left(
\left[  \raisebox{2pt}{\rule{2pt}{1pt}\rule{2pt}{0pt}\rule{2pt}{1pt}\rule
{1pt}{0pt}\raisebox{5pt}{\makebox[0pt]{\hss$Z$\hss}}\raisebox{-11pt}{\makebox
[0pt]{\hss$\tilde{W}$\hss}}\rule{1pt}{0pt}\rule{2pt}{1pt}\rule{2pt}{0pt}%
\rule{2pt}{1pt}}\right]  \right)  =f\left(  \left[  \raisebox{2pt}{\rule
{2pt}{1pt}\rule{2pt}{0pt}\rule{2pt}{1pt}\rule{2pt}{0pt}\rule{2pt}{1pt}%
\rule{2pt}{0pt}\rule{2pt}{1pt}\rule{2pt}{0pt}\rule{2pt}{1pt}\rule{2pt}%
{0pt}\rule{2pt}{1pt}\rule{1pt}{0pt}\raisebox{5pt}{\makebox[0pt]{\hss
$Zg_{1}$\hss}}\raisebox{-11pt}{\makebox[0pt]{\hss$\tilde{W}s_{k}%
g_{2}^{\checkmark}s_{k}$\hss}}\rule{1pt}{0pt}\rule{2pt}{1pt}\rule{2pt}%
{0pt}\rule{2pt}{1pt}\rule{2pt}{0pt}\rule{2pt}{1pt}\rule{2pt}{0pt}\rule
{2pt}{1pt}\rule{2pt}{0pt}\rule{2pt}{1pt}\rule{2pt}{0pt}\rule{2pt}{1pt}%
}\right]  \right)  , \label{eqFD.47}%
\end{equation}
where $Z\in\mathbb{C}^{p\times k}$, $\tilde{W}\in\mathbb{C}^{q\times k}$,
$p+q=n$, $g_{1},g_{2}\in\mathrm{U}\left(  k\right)  $. Then $G^{\prime
}=\mathrm{U}\left(  p\right)  \times\mathrm{U}\left(  q\right)  $ acts on
$\mathcal{F}_{n\times k}$ via the outer tensor product%
\begin{equation}
\left[  R_{\mathrm{U}\left(  p\right)  }^{\prime}\mathop{\hat{\otimes}%
}R_{\mathrm{U}\left(  q\right)  ^{\checkmark}}^{\prime}\right]  \left(
g_{1}^{\prime},g_{2}^{\prime}\right)  f\left(  \left[  \raisebox{2pt}%
{\rule{2pt}{1pt}\rule{2pt}{0pt}\rule{2pt}{1pt}\rule{1pt}{0pt}\raisebox
{5pt}{\makebox
[0pt]{\hss$Z$\hss}}\raisebox{-11pt}{\makebox[0pt]{\hss$\tilde{W}$\hss}}%
\rule{1pt}{0pt}\rule{2pt}{1pt}\rule{2pt}{0pt}\rule{2pt}{1pt}}\right]  \right)
=f\left(  \left[  \raisebox{2pt}{\rule{2pt}{1pt}\rule{2pt}{0pt}\rule{2pt}%
{1pt}\rule{2pt}{0pt}\rule{2pt}{1pt}\rule{2pt}{0pt}\rule{2pt}{1pt}\rule
{2pt}{0pt}\rule{2pt}{1pt}\rule{2pt}{0pt}\rule{2pt}{1pt}\rule{2pt}{0pt}%
\rule{2pt}{1pt}\rule{2pt}{0pt}\rule{2pt}{1pt}\rule{1pt}{0pt}\raisebox
{5pt}{\makebox[0pt]{\hss$\left( g_{1}^{\prime}\right) ^{t}Z$\hss}}%
\raisebox{-11pt}{\makebox[0pt]{\hss$s_{q}\left( g_{2}^{\prime}\right
) ^{-1}s_{q}\tilde{W}$\hss}}\rule{1pt}{0pt}\rule{2pt}{1pt}\rule{2pt}{0pt}%
\rule{2pt}{1pt}\rule{2pt}{0pt}\rule{2pt}{1pt}\rule{2pt}{0pt}\rule{2pt}%
{1pt}\rule{2pt}{0pt}\rule{2pt}{1pt}\rule{2pt}{0pt}\rule{2pt}{1pt}\rule
{2pt}{0pt}\rule{2pt}{1pt}\rule{2pt}{0pt}\rule{2pt}{1pt}}\right]  \right)  ,
\label{eqFD.48}%
\end{equation}
where $\left(  g_{1}^{\prime},g_{2}^{\prime}\right)  \in\mathrm{U}\left(
p\right)  \times\mathrm{U}\left(  q\right)  $.

It follows that we have the isotypic decomposition for the dual pairs $\left(
G^{\prime},G\right)  $%
\begin{equation}
\mathcal{F}_{n\times k}=\sum\limits_{\left(  \nu\right)  \otimes\left(
\lambda^{\checkmark}\right)  }{\!{\oplus}\,}\mathcal{I}_{n\times k}^{\left(
\nu\right)  \otimes\left(  \lambda^{\checkmark}\right)  }, \label{eqFD.49}%
\end{equation}
where $\mathcal{I}_{n\times k}^{\left(  \nu\right)  \otimes\left(
\lambda^{\checkmark}\right)  }$ is isomorphic to $\mathcal{I}_{p\times
k}^{\left(  \nu\right)  }\otimes\mathcal{I}_{q\times k}^{\left(
\lambda^{\checkmark}\right)  }$.

Let $H=\left\{  \left(  g,g\right)  :g\in\mathrm{U}\left(  k\right)  \right\}
$; then $H$ is isomorphic to $\mathrm{U}\left(  k\right)  $ and $H$ acts on
$\mathcal{F}_{n\times q}$ via the inner \textup{(}or Kronecker\textup{)}
tensor product $R_{H}=R_{\mathrm{U}\left(  k\right)  }\otimes R_{\mathrm{U}%
\left(  k\right)  ^{\checkmark}}$. Let $J_{n\times k}$ denote the ring of all
$H$ \textup{(}or $H_{\mathbb{C}}\approx\mathrm{GL}_{k}\left(  \mathbb{C}%
\right)  $\textup{)}-invariant polynomials in $\mathcal{P}_{n\times k}$. Then
from \cite{Ton76b} and \cite{Ton95} $J_{n\times k}$ is generated by the
constants and the $p\times q$ algebraically independent polynomials%
\begin{equation}
p_{\alpha\beta}\left(  \left[  \raisebox{2pt}{\rule{2pt}{1pt}\rule{2pt}%
{0pt}\rule{2pt}{1pt}\rule{1pt}{0pt}\raisebox{5pt}{\makebox[0pt]{\hss$Z$\hss}%
}\raisebox{-11pt}{\makebox[0pt]{\hss$\tilde{W}$\hss}}\rule{1pt}{0pt}\rule
{2pt}{1pt}\rule{2pt}{0pt}\rule{2pt}{1pt}}\right]  \right)  =\left(
Zs_{k}\tilde{W}^{t}\right)  _{\alpha\beta}=\sum_{i=1}^{k}Z_{\alpha i}W_{\beta
i},\qquad1\leq\alpha\leq p,\;1\leq\beta\leq q. \label{eqFD.50}%
\end{equation}
It follows that the ring of all $H$ or $\left(  H_{\mathbb{C}}\right)
$-invariant differential operators with constant coefficients is generated by
the constants and the Laplacians%
\begin{equation}
\bigtriangleup_{\alpha\beta}=p_{\alpha\beta}\left(  D\right)  =\sum_{i=1}%
^{k}\frac{\partial^{2}\;}{\partial Z_{\alpha i}\partial W_{\beta i}}%
,\qquad1\leq\alpha\leq p,\;1\leq\beta\leq q. \label{eqFD.51}%
\end{equation}
Together with the infinitesimal action of $\mathrm{GL}_{n}\left(
\mathbb{C}\right)  $ on $\mathcal{F}_{n\times k}$ the $p_{\alpha\beta}$'s and
$\bigtriangleup_{\alpha\beta}$'s generate a Lie algebra isomorphic to
\textrm{su}$\left(  p,q\right)  $ with commutation relations given by
Eq.\ \textup{(6.4)} in \cite{Ton95}. The global action of this infinitesimal
action defines a representation $R_{H^{\prime}}^{\prime}$ of $H^{\prime
}=\mathrm{SU}\left(  p,q\right)  $ on $\mathcal{F}_{n\times k}$ which is dual
to the representation $R_{H}$.

An element $p$ of $\mathcal{P}_{n\times k}$ is called $H$\emph{-harmonic} if
$\bigtriangleup_{\alpha\beta}p=0$ for all $\alpha=1,\dots,p$, and
$\beta=1,\dots,q$. Let $\mathcal{H}_{n\times k}$ denote the subspace of all
$H$-harmonic polynomial functions of $\mathcal{P}_{n\times k}$ and let
$\mathcal{H}_{n\times k}\left(  \mu\right)  $ denote the subspace of
$\mathcal{H}_{n\times k}$ generated by the elements $f\in\mathcal{P}_{p\times
k}^{\left(  \nu\right)  }\otimes\mathcal{P}_{q\times k}^{\left(
\lambda^{\checkmark}\right)  }$ which also satisfy the condition
$\bigtriangleup_{\alpha\beta}=0$, $1\leq\alpha\leq p$, $1\leq\beta\leq q$. Let
$R_{H}^{\left(  \mu\right)  }$, $\mu=\left(  \nu\right)  \otimes\left(
\lambda^{\checkmark}\right)  $, denote the representation of $H$ on
$\mathcal{H}_{n\times k}\left(  \mu\right)  $ defined by%
\begin{equation}
\left[  R_{H}^{\left(  \mu\right)  }\left(  g\right)  f\right]  \left(
\left[  \raisebox{2pt}{\rule{2pt}{1pt}\rule{2pt}{0pt}\rule{2pt}{1pt}\rule
{1pt}{0pt}\raisebox{5pt}{\makebox[0pt]{\hss$Z$\hss}}\raisebox{-11pt}%
{\makebox[0pt]{\hss$\tilde{W}$\hss}}\rule{1pt}{0pt}\rule{2pt}{1pt}\rule
{2pt}{0pt}\rule{2pt}{1pt}}\right]  \right)  =f\left(  \left[  \raisebox
{2pt}{\rule{2pt}{1pt}\rule{2pt}{0pt}\rule{2pt}{1pt}\rule{2pt}{0pt}\rule
{2pt}{1pt}\rule{2pt}{0pt}\rule{2pt}{1pt}\rule{2pt}{0pt}\rule{2pt}{1pt}%
\rule{2pt}{0pt}\rule{2pt}{1pt}\rule{1pt}{0pt}\raisebox{5pt}{\makebox
[0pt]{\hss$Zg$\hss}}\raisebox{-11pt}{\makebox[0pt]{\hss$\tilde{W}%
s_{k}g^{\checkmark}s_{k}$\hss}}\rule{1pt}{0pt}\rule{2pt}{1pt}\rule{2pt}%
{0pt}\rule{2pt}{1pt}\rule{2pt}{0pt}\rule{2pt}{1pt}\rule{2pt}{0pt}\rule
{2pt}{1pt}\rule{2pt}{0pt}\rule{2pt}{1pt}\rule{2pt}{0pt}\rule{2pt}{1pt}%
}\right]  \right)  , \label{eqFD.52}%
\end{equation}
for all $g\in H$. Then Theorem \textup{5.2} of \cite{Ton95} implies that:

\emph{The representation }$R_{H}^{\left(  \mu\right)  }$\emph{ of }%
$H\approx\mathrm{U}\left(  k\right)  $\emph{ on }$\mathcal{H}_{n\times
k}\left(  \mu\right)  $\emph{ is an irreducible unitary representation of
class }$\left(  \mu\right)  $\emph{ which has signature}%
\begin{equation}
\left(  \mu\right)  =(\underset{k}{\underbrace{\nu_{1},\dots,\nu_{p}%
,0,\dots,0,-\lambda_{q},\dots,-\lambda_{1}}}), \label{eqFD.53}%
\end{equation}
where in Eq.\ \textup{(\ref{eqFD.53})} $\nu_{\alpha}$, $1\leq\alpha\leq p$,
and $\lambda_{\beta}$, $1\leq\beta\leq q$, are integers such that $\nu_{1}%
\geq\dots\geq\nu_{p}\geq0$ and $\lambda_{1}\geq\dots\geq\lambda_{q}\geq0$. Let
$f_{\mu}\left(  \left[  \raisebox{2pt}{\rule{2pt}{1pt}\rule{2pt}{0pt}%
\rule{2pt}{1pt}\rule{1pt}{0pt}\raisebox{5pt}{\makebox[0pt]{\hss$Z$\hss}%
}\raisebox{-11pt}{\makebox[0pt]{\hss$\tilde{W}$\hss}}\rule{1pt}{0pt}\rule
{2pt}{1pt}\rule{2pt}{0pt}\rule{2pt}{1pt}}\right]  \right)  =f_{\nu}\left(
Z\right)  f_{\lambda^{\checkmark}}\left(  \tilde{W}\right)  $, where $f_{\nu}$
is given by Eq.\ \textup{(\ref{eqFD.30})} with $\nu$ replacing $\lambda$ and
$f_{\lambda^{\checkmark}}$ is given by Eq.\ \textup{(\ref{eqFD.44}).} Let
$\mathcal{I}_{n\times k}^{\left(  \mu\right)  }$ be the $H^{\prime}\times
H$-cyclic module generated by $f_{\mu}$; then a proof similar to the previous
cases shows that $\mathcal{H}_{n\times k}^{\prime}\left(  \mu\right)
\mathop{\hat{\otimes}}\mathcal{H}_{n\times k}^{{}}\left(  \mu\right)  $ is
isomorphic to $\mathcal{I}_{n\times k}^{\left(  \mu\right)  }$. By the
``separation of variables theorem'' \textup{1.5} of \cite{Ton76b} and Theorem
\textup{5.1} of \cite{Ton95} it follows that the orthogonal direct sum
decomposition $\mathcal{F}_{n\times k}=\sum\limits_{\left(  \mu\right)
}{\!{\oplus}\,}\mathcal{I}_{n\times k}^{\left(  \mu\right)  }$ holds.
Therefore the reciprocity theorem \textup{\ref{ThmFD.2}} also holds for these
pairs $\left(  G^{\prime},G\right)  $ and $\left(  H^{\prime},H\right)  $. 
\medskip

%\item \label{ExaFD.3(4)}
%
%
%
%
%
%
\noindent4) This example is a generalization of the previous example. Consider
$r$ copies of one of the following groups: $\mathrm{U}\left(  k\right)  $,
$\mathrm{SO}\left(  k\right)  $, or \textrm{Sp}$\left(  k\right)  $, with $k$
even for the last, and let each of them act on a Bargmann--Segal--Fock space
$\mathcal{F}_{p_{i}\times k}$, $1\leq i\leq r$, by right translations. Let
$p_{1}+p_{2}+\dots+p_{r}=n$, and let $G$ denote the direct product of $r$
copies of each type of group. In the case of $\mathrm{U}\left(  k\right)  $ we
allow the $r^{\text{th}}$ copy to act on $\mathcal{F}_{p_{2}\times k}$ either
directly or contragrediently; for the other cases it is not necessary to
consider the contragredient representations since they are identical to the
direct representations.

On each $\mathcal{F}_{p_{i}\times k}$ for the $\mathrm{U}\left(  k\right)  $
action we have the dual action of $\mathrm{U}\left(  p_{i}\right)  $ by left
translations, and with possibly the dual \textup{(}left\textup{)}
contragredient representation in the case $i=r$. For $\mathrm{SO}\left(
k\right)  $ we have the metaplectic representation of $\widetilde
{\mathrm{Sp}_{2p_{i}}}\left(  \mathbb{R}\right)  $, and for \textrm{Sp}%
$\left(  k\right)  $ we have the corresponding representation of
$\mathrm{SO}^{\ast}\left(  2p_{i}\right)  $. Let $G^{\prime}$ denote the dual
group of $G$ thus obtained. Let $H$ denote the diagonal subgroup of $G$; then
in the case of $\mathrm{U}\left(  k\right)  $ an element of $H$ is of the form
$(\underset{r}{\underbrace{u,u,\dots,u}})$ or $(\underset{r-1}{\underbrace
{u,\dots,u},\bar{u}})$, $u\in\mathrm{U}\left(  k\right)  $, and in other cases
an element of $H$ is of the form $(\underset{r}{\underbrace{u,u,\dots,u}})$,
$u\in\mathrm{SO}\left(  k\right)  $ or $u\in$\textrm{Sp}$\left(  k\right)  $.
Let $H^{\prime}$ denote the dual group of $H$ thus obtained. Then $H^{\prime}$
is isomorphic in each case to $\mathrm{U}\left(  n\right)  $\textrm{,
}$\widetilde{\mathrm{Sp}_{2n}}\left(  \mathbb{R}\right)  $, or $\mathrm{SO}%
^{\ast}\left(  2n\right)  $. As in previous examples it is straightforward to
verify that the reciprocity theorem \textup{\ref{ThmFD.2}} holds for these
pairs $\left(  G^{\prime},G\right)  $ and $\left(  H^{\prime},H\right)  $.
%\end{enumerate}
%
%
%
%
%
%
\end{examples}

\section{\label{FID}Reciprocity Theorems for Finite-Infinite Dimensional Dual
Pairs of Groups}

Let $\mathcal{H}$ be an infinite-dimensional separable complex Hilbert space
with a fixed basis $\left\{  e_{1},e_{2},\dots,e_{k},\dots\right\}  $. Let
$\mathrm{GL}_{k}\left(  \mathbb{C}\right)  $ denote the group of all
invertible bounded linear operators on $\mathcal{H}$ which leave the vectors
$e_{n}$, $n>k$, fixed. We define $\mathrm{\ GL}_{\infty}\left(  \mathbb{C}%
\right)  $ as the \emph{inductive limit} of the ascending chain of subgroups%
\[
\mathrm{GL}_{1}\left(  \mathbb{C}\right)  \subset\dots\subset\mathrm{GL}%
_{k}\left(  \mathbb{C}\right)  \subset\cdots.
\]
Thus%
\begin{multline*}
\mathrm{GL}_{\infty}\left(  \mathbb{C}\right)  =\{A=\left(  a_{ij}\right)
,\;i,j\in\mathbb{N}\mid A\text{ is invertible }\\
\text{and all but a finite number of }a_{ij}-\delta_{ij}\text{ are }0\}.
\end{multline*}
If for each $k$ we have a Lie subgroup $G_{k}$ of $\mathrm{GL}_{k}\left(
\mathbb{C}\right)  $ such that $G_{k}$ is naturally embedded in $G_{k+1}$,
$k=1,\dots,n,\dots$, then we can define the \emph{inductive limit} $G_{\infty
}=\varinjlim G_{k}=\bigcup_{k=1}^{\infty}G_{k}$. For example, $\mathrm{U}%
\left(  \infty\right)  =\left\{  u\in\mathrm{GL}_{\infty}\left(
\mathbb{C}\right)  :u^{\ast}=u^{-1}\right\}  $, and thus $\mathrm{U}\left(
\infty\right)  $ is the inductive limit of the groups $\mathrm{U}_{k}$ of all
unitary operators of $\mathcal{H}$ which leave the vectors $e_{n}$, $n>k$, 
fixed.

Following Ol'shanskii we call a unitary representation of $G_{\infty}$
\emph{tame} if it is continuous in the group topology in which the ascending
chain of subgroups of type $\left\{
\begin{pmatrix}
1_{k} & 0\\
0 & \ast
\end{pmatrix}
\right\}  ,$ $k=1,2,3,\dots$, constitutes a fundamental system of
neighborhoods of the identity $1_{\infty}$. Assume that for each $k$ a
continuous unitary representation $\left(  R_{k},\mathcal{H}_{k}\right)  $ is
given and an isomorphic embedding $i_{k+1}^{k}\colon\mathcal{H}_{k}%
\rightarrow\mathcal{H}_{k+1}$ commuting with the action of $G_{k}$ (i.e.,
$i_{k+1}^{k}\circ R_{k}\left(  g\right)  =R_{k+1}\left(  g\right)  \circ
i_{k+1}^{k}$) is given. For $j\leq k$ define the \emph{connecting map}
$\varphi_{jk}\colon G_{j}\times\mathcal{H}_{j}\rightarrow G_{k}\times
\mathcal{H}_{k}$ by%
\begin{equation}
\varphi_{jk}\left(  g_{j},x_{j}\right)  =\left(  g_{k},x_{k}\right)
,\qquad\left(  g_{j},x_{j}\right)  \in G_{j}\times\mathcal{H}_{j},
\label{eqFID.1}%
\end{equation}
where in Eq.\ (\ref{eqFID.1}) $g_{k}$ (resp.\ $x_{k}$) denotes the natural
embedding of $g_{j}$ (resp.\ $x_{j}$) in $G_{k}$ (resp.\ $\mathcal{H}_{k}$).
Then obviously the diagram%
\begin{equation}%
\begin{array}
[c]{ccc}%
G_{j}\times\mathcal{H}_{j} & \overset{R_{j}}{\longrightarrow} & \mathcal{H}%
_{j}\\
\smash{\vcenter{\hrule width0.275pt height14pt}\kern-0.275pt\makebox[0.275pt]
{\hss
\raisebox{-3pt}{$\downarrow$}\hss}}\vrule width0ptheight22ptdepth12pt\llap
{$\scriptstyle\varphi_{jk}\;$} &  & \smash{\vcenter{\hrule
width0.275pt height14pt}\kern-0.275pt\makebox[0.275pt]
{\hss\raisebox{-3pt}{$\downarrow$}\hss}%
}\vrule width0ptheight22ptdepth12pt\rlap{$\scriptstyle\;i_{k}^{j}=i_{k}%
^{k-1}\circ\dots\circ i_{j+2}^{j+1}\circ i_{j+1}^{j}$}\\
G_{k}\times\mathcal{H}_{k} & \overset{R_{k}}{\longrightarrow} & \mathcal{H}%
_{k}%
\end{array}
\label{eqFID.2}%
\end{equation}
is commutative. Let $\mathcal{H}_{\infty}$ denote the Hilbert-space completion
of $\bigcup_{k=1}^{\infty}\mathcal{H}_{k}$ and define a representation
$R_{\infty}$ of $G_{\infty}$ on $\mathcal{H}_{k}$ by%
\begin{equation}
R_{\infty}\left(  g\right)  x=R_{k}\left(  g\right)  x\text{\qquad if }g\in
G_{k}\text{ and }x\in\mathcal{H}_{k}. \label{eqFID.3}%
\end{equation}
Then obviously $R_{\infty}$ is a unique continuous unitary representation of
$G_{\infty}$ on $\bigcup_{k=1}^{\infty}\mathcal{H}_{k}$ which can be extended
to a unique continuous unitary representation of $G_{\infty}$ on
$\mathcal{H}_{\infty}$. Let $\varphi_{k}$ denote the canonical map of $\left(
G_{k},\mathcal{H}_{k}\right)  $ into $\left(  G_{\infty},\mathcal{H}_{\infty
}\right)  $ and $i_{k}$ denote the canonical map of $\mathcal{H}_{k}$ into
$\mathcal{H}_{\infty}$; then obviously the diagram
\begin{equation}%
\begin{array}
[c]{ccc}%
G_{k}\times\mathcal{H}_{k} & \overset{R_{k}}{\longrightarrow} & \mathcal{H}%
_{k}\\
\smash{\vcenter{\hrule width0.275pt height14pt}\kern-0.275pt\makebox[0.275pt]
{\hss
\raisebox{-3pt}{$\downarrow$}\hss}}\vrule width0ptheight22ptdepth12pt\llap
{$\scriptstyle\varphi_{k}\;$} &  & \smash{\vcenter{\hrule
width0.275pt height14pt}\kern-0.275pt\makebox[0.275pt]
{\hss\raisebox{-3pt}{$\downarrow$}\hss}%
}\vrule width0ptheight22ptdepth12pt\rlap{$\scriptstyle\;i_{k}$}\\
G_{\infty}\times\mathcal{H}_{\infty} & \overset{R_{\infty}}{\longrightarrow} &
\mathcal{H}_{\infty}%
\end{array}
\label{eqFID.4}%
\end{equation}
is commutative.

The following theorem, which is well-known when $i_{k+1}^{k}$ is an isometric
embedding (see, e.g., \cite{Ol'90}), is crucial for what follows.

\begin{theorem}
\label{ThmFID.1}If the representations $\left(  R_{k},\mathcal{H}_{k}\right)
$ are all irreducible then the inductive limit representation $\left(
R_{\infty},\mathcal{H}_{\infty}\right)  $ is also irreducible.
\end{theorem}

\begin{proof}
Let $A$ be a bounded operator on $\mathcal{H}_{\infty}$ which belongs to the
commutant of the algebra of operators generated by the set $\left\{
R_{\infty}\left(  g\right)  ,g\in G_{\infty}\right\}  $. Since $\bigcup
_{k=1}^{\infty}\mathcal{H}_{k}$ is dense in $\mathcal{H}_{\infty}$ and all the
linear operators involved are continuous we can without loss of generality
consider them as operating on $\bigcup_{k=1}^{\infty}\mathcal{H}_{k}$ and
satisfying $A\left(  i_{l}^{k}\left(  x\right)  \right)  =Ax$ for $k\leq l$
and for all $x\in\mathcal{H}_{k}$. Let $P_{k}$ denote the projection of
$\bigcup_{k=1}^{\infty}\mathcal{H}_{k}$ onto $\mathcal{H}_{k}$. Let $A_{k}$
denote the restriction of $A$ to $\mathcal{H}_{k}$; then $A_{k}$ is a bounded
linear operator of $\mathcal{H}_{k}$ into $\bigcup_{n=1}^{\infty}%
\mathcal{H}_{n}$. It follows immediately that $P_{k}A_{k}\colon\mathcal{H}%
_{k}\rightarrow\mathcal{H}_{k}$ is a bounded linear operator on $\mathcal{H}%
_{k}$. Let $x\in\mathcal{H}_{k}$ and suppose $A_{k}x=Ax$ belongs to
$\mathcal{H}_{l}$. If $l\leq k$ we may use the isomorphic embedding $i_{k}%
^{l}=i_{k}^{k-1}\circ\dots\circ i_{l+1}^{l}\colon\mathcal{H}_{l}%
\rightarrow\mathcal{H}_{k}$ to identify $Ax$ with an element of $\mathcal{H}%
_{k}$ so that $P_{k}A_{k}x=A_{k}x=Ax$, and thus%
\[
R_{k}\left(  g_{k}\right)  P_{k}A_{k}x=R_{\infty}\left(  g_{k}\right)
Ax=AR_{\infty}\left(  g_{k}\right)  x=P_{k}R_{k}\left(  g_{k}\right)
x,\qquad\forall\,g_{k}\in G_{k}.
\]
If $l>k$ then use $i_{l}^{k}$ to identify $\mathcal{H}_{k}$ with a subspace of
$\mathcal{H}_{l}$. Write $Ax=y+z$ where $y$ belongs to the identified subspace
of $\mathcal{H}_{k}$ and $z$ belongs to its orthogonal complement in
$\mathcal{H}_{l}$. Since all representations are unitary and for $g_{k}\in
G_{k}$ we have $i_{l}^{k}\circ R_{k}\left(  g_{k}\right)  =R_{l}\left(
g_{k}\right)  \circ i_{l}^{k}$ it follows that%
\[
P_{k}R_{\infty}\left(  g_{k}\right)  A_{k}x=P_{k}R_{k}\left(  g_{k}\right)
y=R_{\infty}\left(  g_{k}\right)  P_{k}Ax.
\]
By assumption $R_{\infty}\left(  g_{k}\right)  Ax=AR_{\infty}\left(
g_{k}\right)  x$, therefore%
\begin{multline*}
R_{k}\left(  g_{k}\right)  P_{k}A_{k}x=R_{\infty}\left(  g_{k}\right)
P_{k}Ax\\
=P_{k}R_{\infty}\left(  g_{k}\right)  A_{k}x=P_{k}A_{k}R_{\infty}\left(
g_{k}\right)  x=P_{k}A_{k}R_{k}\left(  g_{k}\right)  x.
\end{multline*}
Since this relation holds for all $x\in\mathcal{H}_{k}$ and $g_{k}\in G_{k}$
it follows that $P_{k}A_{k}$ belongs to the commutant of the algebra of
operators on $\mathcal{H}_{k}$ generated by the set $\left\{  R_{k}\left(
g_{k}\right)  ,g_{k}\in G_{k}\right\}  $. Schur's lemma for operator algebras
(see, e.g., \cite[Proposition 2.3.1, p.~39]{Dix69}) implies that $P_{k}%
A_{k}=\lambda_{k}I_{k}$, where $\lambda_{k}$ is a scalar depending on $k$ and
$I_{k}$ is the identity operator on $\mathcal{H}_{k}$. Now $A$ is a map of
inductive limit sets such that $P_{k}A_{k}\colon\mathcal{H}_{k}\rightarrow
\mathcal{H}_{k}$, and it follows from the definition of an inductive limit map
that $\lambda_{k}=\lambda_{l}$ for sufficiently large $k$, $l$ with $k<l$.
Indeed, if $x\in\mathcal{H}_{k}$ and $A_{k}x=Ax\in\mathcal{H}_{j}$ with $j\leq
k$ then $P_{k}A_{k}x=i_{k}^{j}\left(  Ax\right)  =\lambda_{k}x$. For $l>k$ we
then have%
\begin{multline*}
\lambda_{l}^{{}}i_{l}^{k}\left(  x\right)  =P_{l}A_{l}\left(  i_{l}^{k}\left(
x\right)  \right)  =P_{l}A\left(  i_{l}^{k}\left(  x\right)  \right)
=P_{l}Ax\\
=i_{l}^{j}\left(  Ax\right)  =i_{l}^{k}\left(  i_{k}^{j}\left(  Ax\right)
\right)  =i_{l}^{k}\left(  P_{k}A_{k}\left(  x\right)  \right)  =\lambda
_{k}^{{}}i_{l}^{k}\left(  x\right)  .
\end{multline*}
On the other hand, if $Ax\in\mathcal{H}_{j}$ with $j>k$ then for all $l\geq j$
we have%
\begin{multline*}
P_{l}A_{l}\left(  i_{l}^{k}\left(  x\right)  \right)  =P_{l}A\left(  i_{l}%
^{k}\left(  x\right)  \right)  =P_{l}Ax=P_{l}A_{j}\left(  i_{j}^{k}\left(
x\right)  \right) \\
=P_{l}P_{j}A_{j}\left(  i_{j}^{k}\left(  x\right)  \right)  =P_{l}\left(
\lambda_{j}^{{}}i_{j}^{k}\left(  x\right)  \right)  =i_{l}^{j}\left(
\lambda_{j}^{{}}i_{j}^{k}\left(  x\right)  \right)  =\lambda_{j}^{{}}i_{l}%
^{j}\left(  i_{j}^{k}\left(  x\right)  \right)  =\lambda_{j}^{{}}i_{l}%
^{k}\left(  x\right)  .
\end{multline*}
Since $P_{l}A_{l}\left(  i_{l}^{k}\left(  x\right)  \right)  =\lambda_{l}^{{}%
}i_{l}^{k}\left(  x\right)  $, we must have $\lambda_{j}=\lambda_{l}$ for all
$l\geq j$. This implies that $A=\lambda I_{\infty}$ where $\lambda
\in\mathbb{C}$ is a constant and $I_{\infty}$ is the identity on
$\mathcal{H}_{\infty}$. By the same Schur's lemma quoted above the
representation $R_{\infty}$ on $\mathcal{H}_{\infty}$ must be irreducible.
\end{proof}

Now fix $n$ and consider the chain of Hilbert spaces $\mathcal{F}_{n\times k}$
from Section \ref{FD} with $k>2n$. Let $\left(  G_{n}^{\prime},G_{k}^{{}%
}\right)  $ denote a dual pair of groups with dual representations $\left(
R_{n}^{\prime},R_{k}^{{}}\right)  $ acting on $\mathcal{F}_{n\times k}$ as in
Theorem \ref{ThmFD.2}. Then we have the chain of embedded subgroups
$G_{k}\subset G_{k+1}\subset\cdots$; for example, $\mathrm{U}\left(  k\right)
$ is naturally embedded in $\mathrm{U}\left(  k+1\right)  $ via the embedding
$u\rightarrow\left(
\begin{smallmatrix}
u & 0\\
0 & 1
\end{smallmatrix}
\right)  $, $u\in\mathrm{U}\left(  k\right)  $. Therefore we can define the
inductive limit $G_{\infty}=\varinjlim G_{k}=\bigcup_{k>2n}^{\infty}G_{k}$. We
also have an isometric embedding $i_{k+1}^{k}\colon\mathcal{F}_{n\times
k}\rightarrow\mathcal{F}_{n\times\left(  k+1\right)  }$ such that%
\[
i_{k+1}^{k}\circ R_{k}\left(  g\right)  =R_{k+1}\left(  g\right)  \circ
i_{k+1}^{k}.
\]
To see this we take the case $n=1$: then an element $f$ of $\mathcal{F}%
_{n\times k}$ is a function of $Z=\left(  Z_{1},\dots,Z_{k}\right)  $ of the
form given by Eq.\ (\ref{eqFD.7}), and the verification of the equation above
is straightforward. Let $\mathcal{F}_{n\times\infty}$ denote the Hilbert-space
completion on $\bigcup_{k>2n}^{\infty}\mathcal{F}_{n\times k}$. Then it is
clear that the inductive limit representation $R_{\infty}$ of $G_{\infty}$ on
$\mathcal{F}_{n\times\infty}$ is tame and satisfies the relations
(\ref{eqFID.2}), (\ref{eqFID.3}), and (\ref{eqFID.4}).

If $G_{k}$ is a compact group then every irreducible unitary representation of
$G_{k}$ is of the form $\left(  \rho_{\lambda_{k}},V_{\lambda_{k}}\right)  $
with highest weight $\left(  \lambda_{k}\right)  =\left(  m_{1},m_{2}%
,\dots,m_{i},\dots\right)  $, where $m_{1},m_{2},\dots$ are nonnegative
integers satisfying $m_{1}\geq m_{2}\geq\cdots$ and the numbers $m_{i}$ are
equal to $0$ for sufficiently large $i$. Consider the decomposition
(\ref{eqFD.5}) of Definition \ref{DefFD.1} of the dual module $\mathcal{F}%
_{n\times k}$ into isotypic components%
\[
\mathcal{F}_{n\times k}=\sum\limits_{\makebox[0pt]{\hss$\scriptstyle
\left(
\lambda_{k}\right)  $\hss}}{\!{\oplus}\,}\mathcal{I}_{n\times k}^{\left(
\lambda_{k}\right)  }%
\]
where the signatures $\left(  \lambda_{k}\right)  $ actually depend
essentially on $n$, but since $n$ is fixed, to alleviate the notation we just
tacitly assume this dependence. Also for $k$ sufficiently large if $\left(
\lambda_{k}\right)  =\left(  m_{1},\dots,m_{i},\dots\right)  $ then $\left(
\lambda_{k+1}\right)  =\left(  m_{1},\dots,m_{i},\dots,\dots\right)  $ and we
write succinctly $\left(  \lambda_{k}\right)  \subset\left(  \lambda
_{k+1}\right)  $.

For sufficiently large $k$ we can exhibit an isomorphic embedding $i_{k+1}%
^{k}\colon\mathcal{I}_{n\times k}^{\left(  \lambda_{k}\right)  }%
\rightarrow\mathcal{I}_{n\times\left(  k+1\right)  }^{\left(  \lambda
_{k+1}\right)  }$. If $H_{k}$ is a subgroup of $G_{k}$ such that
$H_{n}^{\prime}$ contains $G_{n}^{\prime}$ and $\left(  H_{n}^{\prime}%
,H_{k}^{{}}\right)  $ forms a dual pair then the same process can be repeated
for the chain $\left(  H_{n}^{\prime},H_{k}^{{}}\right)  \subset\left(
H_{n}^{\prime},H_{k+1}^{{}}\right)  \subset\cdots$. If $G_{k}$ (or $H_{k}$) is
of the type $\underset{r}{\underbrace{\mathrm{U}\left(  k\right)  \times
\dots\times\mathrm{U}\left(  k\right)  }}$ then each $i_{k+1}^{k}$ is an
isometric embedding; for other types of $G_{k}$ (or $H_{k}$) the definition of
$i_{k+1}^{k}$ is more subtle. This can be examined case by case although the
process is very tedious. To illustrate this we consider the case
$\mathcal{F}_{1\times k}$ with $H_{k}=\mathrm{SO}\left(  k\right)  $ and
$G_{1}=\widetilde{\mathrm{Sp}_{2}\left(  \mathbb{R}\right)  }=\widetilde
{\mathrm{SL}_{2}\left(  \mathbb{R}\right)  }$. Then Eq.\ (\ref{eqFD.22}) and
Eq.\ (\ref{eqFD.23}) imply that%
\[
\mathcal{F}_{1\times k}=\sum\limits_{\makebox[0pt]{\hss$\scriptstyle
r=0$\hss}}^{\infty}{\!{\oplus}\,}\mathcal{I}_{1\times k}^{\left(  r\right)
_{k}}\text{\qquad with\qquad}\mathcal{I}_{1\times k}^{\left(  r\right)  _{k}}=
\sum\limits_{\makebox[0pt]{\hss$\scriptstyle
j=0$\hss}}^{\infty}{\!{\oplus}\,}p_{0,k}^{j}\mathcal{H}_{1\times k}^{\left(
r\right)  _{k}}%
\]
where $p_{0,k}\left(  Z\right)  =Z_{1}^{2}+\dots+Z_{k}^{2}$, $\left(
r\right)  _{k}=(\underset{k}{\underbrace{r,0,\dots,0}})$, and $\mathcal{H}%
_{1\times k}^{\left(  r\right)  }$ are the subspace of all harmonic
homogeneous polynomials of degree $r$. Obviously a harmonic homogeneous
polynomial $h$ of degree $r$ in $k$ variables can be considered as a harmonic
homogeneous polynomial of $r$ in $k+1$ variables. So we can define an
isomorphic embedding $i_{k+1}^{k}\colon\mathcal{I}_{1\times k}^{\left(
r\right)  _{k}}\rightarrow\mathcal{I}_{1\times\left(  k+1\right)  }^{\left(
r\right)  _{k+1}}$ by sending $p_{0,k}^{j}h$ into $p_{0,\left(  k+1\right)
}^{j}h$, and clearly
\begin{multline*}
R_{H}\left(  u_{k}\right)  \left(  p_{0,\left(  k+1\right)  }^{j}h\right)
=p_{0,\left(  k+1\right)  }^{j}R_{H}\left(  u_{k}\right)  h=i_{k+1}^{k}%
p_{0,k}^{j}R_{H}\left(  u_{k}\right)  h\\
=i_{k+1}^{k}\left(  \left(  R_{H}\left(  u_{k}\right)  p_{0,k}^{j}\right)
\left(  R_{H}\left(  u_{k}\right)  h\right)  \right)  =i_{k+1}^{k}\left(
R_{H}\left(  u_{k}\right)  \left(  p_{0,k}^{j}h\right)  \right)
\end{multline*}
for all $u_{k}\in H_{k}$. Thus, $R_{H}\left(  u_{k}\right)  \circ i_{k+1}%
^{k}=i_{k+1}^{k}\circ R_{H}\left(  u_{k}\right)  $ for all $u_{k}\in H_{k}$.
It follows that $i_{k+1}^{k}$ can be extended to the whole space
$\mathcal{F}_{1\times k}$ and that $i_{k+1}^{k}\left(  \mathcal{F}_{1\times
k}\right)  =\sum\limits_{\makebox[0pt]{\hss$\scriptstyle
r=0$\hss}}^{\infty}{\!{\oplus}\,}i_{k+1}^{k}\left(  \mathcal{I}_{1\times
k}^{\left(  r\right)  _{k}}\right)  $ is an isomorphic embedding of
$\mathcal{F}_{1\times k}$ into $\mathcal{F}_{1\times\left(  k+1\right)  }$.
Also note in this very special case $\left(  r\right)  _{k}\subset\left(
r\right)  _{k+1}$ for all $k>2$ and that no other signatures $\left(
r\right)  _{k+1}$ occur in $\mathcal{F}_{1\times\left(  k+1\right)  }$ without
$\left(  r\right)  _{k}$ occurring in $\mathcal{F}_{1\times k}$; this fact is
an exception and almost never happens in the general case (e.g., $n\geq2$). By
Theorem \ref{ThmFID.1} the tensor product representations $R_{G_{n}^{\prime}%
}^{\left(  \lambda^{\prime}\right)  }\otimes R_{G_{\infty}}^{\left(
\lambda\right)  }$ and $R_{H_{n}^{\prime}}^{\left(  \mu^{\prime}\right)
}\otimes R_{H_{\infty}}^{\left(  \mu\right)  }$ of $G_{n}^{\prime}\times
G_{\infty}^{{}}$ and $H_{n}^{\prime}\times H_{\infty}^{{}}$ on $\mathcal{I}%
_{n\times\infty}^{\left(  \lambda\right)  }$ and $\mathcal{I}_{n\times\infty
}^{\left(  \mu\right)  }$, respectively, are irreducible with signature
$\left(  \lambda\right)  _{\infty}$ and $\left(  \mu\right)  _{\infty}$,
respectively, where if $\left(  \lambda_{k}\right)  =\left(  m_{1},m_{2}%
,\dots,m_{i},\dots\right)  $ then $\left(  \lambda\right)  _{\infty
}=(\underset{\infty}{\underbrace{m_{1},m_{2},\dots,m_{i},\dots,0}},\dots,0)$
and similarly for $\left(  \mu\right)  _{\infty}$. Note that as $n$ is fixed,
the group $G_{n}^{\prime}$ remains fixed; however, its representation
$R_{G_{n}^{\prime}}^{\prime}$ on $\mathcal{F}_{n\times k}$ does depend on $k$,
and should be written as $\left(  R_{G_{n}^{\prime}}^{\prime}\right)  _{k}$,
and as $k\rightarrow\infty$, $\left(  R_{G_{n}^{\prime}}^{\prime}\right)
_{\infty}$ has to be considered as an inductive limit of representations,
although for $k$ sufficiently large all the representations $\left(
R_{G_{n}^{\prime}}^{\left(  \lambda^{\prime}\right)  }\right)  _{k}$ are
equivalent. The same observations apply to $\left(  R_{H_{n}^{\prime}}%
^{\prime}\right)  _{k}$ and $\left(  R_{H_{n}^{\prime}}^{\left(  \mu^{\prime
}\right)  }\right)  _{k}$. To illustrate this let us consider again the case
$\mathrm{U}\left(  1\right)  \times\mathrm{U}\left(  k\right)  $ and
$\widetilde{\mathrm{SL}_{2}\left(  \mathbb{R}\right)  }\times\mathrm{SO}%
\left(  k\right)  $ acting on $\mathcal{F}_{1\times k}$. Indeed, the
infinitesimal action of $R_{G_{\mathbb{C}}^{\prime}}^{\prime}$ is given by
Eq.\ (\ref{eqFD.11}) as $L_{k}=\sum_{i=1}^{k}Z_{i}\partial/\partial Z_{i}$ and
$L_{k+1}=\sum_{i=1}^{k+1}Z_{i}\partial/\partial Z_{i}$, and for $p\in
\mathcal{P}_{1\times k}^{\left(  m\right)  }\subset\mathcal{P}_{1\times\left(
k+1\right)  }^{\left(  m\right)  }$ Eq.\ (\ref{eqFD.16}) implies that%
\[
L_{k}p=L_{k+1}p=mp.
\]
By Eq.\ (\ref{eqFD.18}) the infinitesimal actions of $R_{H_{1}^{\prime}%
}^{\prime}$ on $\mathcal{F}_{1\times k}$ and $\mathcal{F}_{1\times\left(
k+1\right)  }$ are given, respectively, by
\begin{equation}
\left\{
%TCIMACRO{\TeXButton{aligned equations}{\begin{aligned}
%E_{k} &= \frac{k}{2}+L_{k},
%&\quad X_{k}^{+} &= \frac{1}{2}\sum_{i=1}^{k}Z_{i}^{2},
%&\quad X_{k}^{-} &= \frac{1}{2}\sum_{i=1}^{k}\frac{\partial^{2}\;}{\partial
%Z_{i}^{2}},
%&\quad&\text{and}  \\
%E_{k+1} &= \frac{k+1}{2}+L_{k+1},
%&\quad X_{k+1}^{+} &= \frac{1}{2}\sum_{i=1}^{k+1}Z_{i}^{2},
%&\quad X_{k+1}^{-} &= \frac{1}{2}\sum_{i=1}^{k+1}\frac{\partial^{2}%
%\;}{\partial
%Z_{i}^{2}}.
%&  &
%\end{aligned}
%}}%
%BeginExpansion
\begin{aligned}
E_{k} &= \frac{k}{2}+L_{k},
&\quad X_{k}^{+} &= \frac{1}{2}\sum_{i=1}^{k}Z_{i}^{2},
&\quad X_{k}^{-} &= \frac{1}{2}\sum_{i=1}^{k}\frac{\partial^{2}\;}{\partial
Z_{i}^{2}},
&\quad&\text{and}  \\
E_{k+1} &= \frac{k+1}{2}+L_{k+1},
&\quad X_{k+1}^{+} &= \frac{1}{2}\sum_{i=1}^{k+1}Z_{i}^{2},
&\quad X_{k+1}^{-} &= \frac{1}{2}\sum_{i=1}^{k+1}\frac{\partial^{2}%
\;}{\partial
Z_{i}^{2}}.
&  &
\end{aligned}
%EndExpansion
\right.  \label{eqFID.5}%
\end{equation}
If $h_{k}\in\mathcal{H}_{1\times k}^{\left(  r\right)  }$ then
Eqs.\ (\ref{eqFD.24}), (\ref{eqFD.25}) applied to $\left\{  E_{k}^{{}}%
,X_{k}^{+},X_{k}^{-}\right\}  $ show that $J_{k}h_{k}$ is an irreducible
representation of $\mathrm{sl}_{2}\left(  \mathbb{R}\right)  $ with signature
$\left(  r\right)  $. Similarly if $h_{k+1}\in\mathcal{H}_{1\times\left(
k+1\right)  }^{\left(  r\right)  }$ then $J_{k+1}h_{k+1}$ is also an
irreducible representation of $\mathrm{sl}_{2}\left(  \mathbb{R}\right)  $
with signature $\left(  r\right)  $.

Let $\mathcal{F}_{n\times\infty}$ denote the Hilbert-space completion of
$\bigcup_{k}\mathcal{F}_{n\times k}$; then $\mathcal{F}_{n\times\infty
}=\varinjlim\mathcal{F}_{n\times k}$ is the inductive limit of the chain
$\left\{  \mathcal{F}_{n\times k}\right\}  $.

After this necessary preparatory work we can now state and prove the main
theorem of this paper.

\begin{theorem}
\label{ThmFID.2}Let $G_{\infty}$ denote the inductive limit of a chain
$G_{k}\subset G_{k+1}\subset\cdots$ of compact groups. Let $R_{G_{\infty}}$
and $R_{\left(  G_{n}^{\prime}\right)  _{\infty}}^{\prime}$ be given dual
representations on $\mathcal{F}_{n\times\infty}$. Let $H_{\infty}$ denote the
inductive limit of a chain of compact subgroups $H_{k}\subset H_{k+1}%
\subset\cdots$ such that $H_{k}\subset G_{k}$ for all $k$. Let $R_{H_{\infty}%
}$ be the representation of $H_{\infty}$ on $\mathcal{F}_{n\times\infty}$
obtained by restricting $R_{G_{\infty}}$ to $H_{\infty}$. If there exists a
group $H_{n}^{\prime}\supset G_{n}^{\prime}$ and a representation $R_{\left(
H_{n}^{\prime}\right)  _{\infty}}^{\prime}$ on $\mathcal{F}_{n\times\infty}$
such that $R_{\left(  H_{n}^{\prime}\right)  _{\infty}}^{\prime}$ is dual to
$R_{H_{\infty}}$ and $R_{\left(  G_{n}^{\prime}\right)  _{\infty}}^{\prime}$
is the restriction of $R_{\left(  H_{n}^{\prime}\right)  _{\infty}}^{\prime}$
to the subgroup $G_{n}^{\prime}$ of $H_{n}^{\prime}$ then we have the
following multiplicity-free decompositions of $\mathcal{F}_{n\times\infty}$
into isotypic components:%
\begin{equation}
\mathcal{F}_{n\times\infty}=\sum\limits_{\left(  \lambda\right)  }{\!{\oplus
}\,}\mathcal{I}_{n\times\infty}^{\left(  \lambda\right)  }=\sum
\limits_{\left(  \mu\right)  }{\!{\oplus}\,}\mathcal{I}_{n\times\infty
}^{\left(  \mu\right)  } \label{eqFID.6}%
\end{equation}
where $\left(  \lambda\right)  $ is a common irreducible signature of the pair
$\left(  G_{n}^{\prime},G_{\infty}^{{}}\right)  $ and $\left(  \mu\right)  $
is a common irreducible signature of the pair $\left(  H_{n}^{\prime
},H_{\infty}^{{}}\right)  $.

If $\lambda_{G_{\infty}}$ \textup{(}resp.\ $\lambda_{\left(  G_{n}^{\prime
}\right)  _{\infty}}^{\prime}$\textup{)} denotes an irreducible unitary
representation of class $\left(  \lambda\right)  $ and $\mu^{}_{H_{\infty}}$
\textup{(}resp.\ $\mu_{\left(  H_{n}^{\prime}\right)  _{\infty}}^{\prime}%
$\textup{)} denotes an irreducible unitary representation of class $\left(
\mu\right)  $ then the multiplicity $\dim\left[  \operatorname*{Hom}%
_{H_{\infty}}\left(  \mu^{}_{H_{\infty}}:\lambda_{G_{\infty}}|_{H_{\infty}%
}\right)  \right]  $ of the irreducible representation $\mu^{}_{H_{\infty}}$
in the restriction to $H_{\infty}$ of the representation $\lambda_{G_{\infty}%
}$ is equal to the multiplicity $\dim\left[  \operatorname*{Hom}%
_{G_{n}^{\prime}}\left(  \lambda_{\left(  G_{n}^{\prime}\right)  _{\infty}%
}^{\prime}:\mu_{\left(  H_{n}^{\prime}\right)  _{\infty}}^{\prime}\bigg
|_{G_{n}^{\prime}}\right)  \right]  $ of the irreducible representation
$\lambda_{\left(  G_{n}^{\prime}\right)  _{\infty}}^{\prime}$ in the
restriction to $G_{n}^{\prime}$ of the representation $\mu_{\left(
H_{n}^{\prime}\right)  _{\infty}}^{\prime}$.
\end{theorem}

\begin{proof}
As remarked above, the dual $\left(  G_{n}^{\prime},G_{\infty}^{{}}\right)
$-module $\mathcal{I}_{n\times\infty}^{\left(  \lambda\right)  }$ is
irreducible (by Theorem \ref{ThmFID.1}) with signature $\left(  \lambda
\right)  $, and isotypic components of different signatures are mutually
orthogonal since their projections $\mathcal{I}_{n\times k}^{\left(
\lambda\right)  _{k}}$ are mutually orthogonal. Finally if a vector in
$\mathcal{F}_{n\times\infty}$, which we may assume to belong to $\mathcal{F}%
_{n\times k}$ for some $k$, is orthogonal to $\mathcal{I}_{n\times\infty
}^{\left(  \lambda\right)  }$ for all $\left(  \lambda\right)  $, it must
therefore be orthogonal to $\mathcal{I}_{n\times k}^{\left(  \lambda\right)
_{k}}$ for all $\left(  \lambda\right)  _{k}$, and hence must be the zero
vector in $\mathcal{F}_{n\times k}$, and thus zero in $\mathcal{F}%
_{n\times\infty}$. A similar argument applies to the isotypic components
$\mathcal{I}_{n\times\infty}^{\left(  \mu\right)  }$, and thus
Eq.\ (\ref{eqFID.6}) holds.

Now fix $\left(  \lambda\right)  $ and $\left(  \mu\right)  $. Then the
restriction of $R_{G_{\infty}}$ to $\mathcal{I}_{n\times\infty}^{\left(
\lambda\right)  }$ decomposes into a (non-canonical) orthogonal direct sum of
equivalent irreducible unitary representations of signature $\left(
\lambda\right)  _{\infty}$. A representative of this representation may be
obtained by applying Theorem \ref{ThmFID.1} to get the inductive limit
$\left(  G_{\infty},R_{\left(  \lambda\right)  _{\infty}}\right)  $ of the
chain $\left(  G_{k},R_{\lambda_{k}}\right)  $; for example, when
$G_{k}=\mathrm{U}\left(  k\right)  $, the representation $R_{\lambda_{k}}$ is
given by Eq.\ (\ref{eqFD.29}) on $\mathcal{P}_{n\times k}^{\left(
\lambda\right)  _{k}}$. Considered as a $G_{n}^{\prime}$-module $\mathcal{I}%
_{n\times\infty}^{\left(  \lambda\right)  }$ decomposes into a (non-canonical)
orthogonal direct sum of equivalent irreducible unitary representations of
signature $\left(  \lambda^{\prime}\right)  _{n}$. A representative of this
representation may be obtained by applying Theorem \ref{ThmFID.1} to get the
inductive limit $\left(  G_{n}^{\prime},R_{\left(  \lambda_{n}^{\prime
}\right)  _{\infty}}^{\prime}\right)  $ (note that although $G_{n}^{\prime}$
is a stationary chain at $n$, the representations $R_{\left(  \lambda
_{n}^{\prime}\right)  _{k}}^{\prime}$ depend on $k$ even though they are all
equivalent and belong to the class $\left(  \lambda^{\prime}\right)  _{n}$);
for example, when $G_{n}=\mathrm{U}\left(  n\right)  $ the representation
$R_{\lambda_{n}^{\prime}}^{\prime}$ is given by Eq.\ (\ref{eqFD.32}) on
$\mathcal{P}_{n\times k}^{\left(  \lambda^{\prime}\right)  _{n}}$ which is
defined by Eq.\ (\ref{eqFD.31}). By an analogous argument we infer that the
same conclusions hold for $\left(  \mu\right)  $, $\mathcal{I}_{n\times\infty
}^{\left(  \mu\right)  }$, $\left(  H_{\infty},R_{\left(  \mu\right)
_{\infty}}\right)  $, $\left(  H_{n}^{\prime},R_{\left(  \mu_{n}^{\prime
}\right)  _{\infty}}^{\prime}\right)  $.

Now consider the decomposition of the restriction to $H_{k}$ of the
representation $R_{\lambda_{k}}$ of $G_{k}$. The multiplicity of $\left(
\mu\right)  _{k}$ in $\left(  \lambda\right)  _{k}|_{H_{k}}$ is the dimension
of \linebreak $\operatorname*{Hom}_{H_{k}}\left(  R_{\mu^{}_{k}}%
:R_{\lambda_{k}}|_{H_{k}}\right)  $, where $\operatorname*{Hom}_{H_{k}}\left(
R_{\mu^{}_{k}}:R_{\lambda_{k}}|_{H_{k}}\right)  $ is the vector space of
linear homomorphisms intertwining $R_{\mu^{}_{k}}$ and $R_{\lambda_{k}%
}|_{H_{k}}$. Since $G_{k}$ and $H_{k}$ are, by assumption, compact, this
dimension is finite. If $T_{k}\colon\mathcal{H}_{\mu^{}_{k}}\rightarrow
\mathcal{H}_{\lambda_{k}}$ is an element of $\operatorname*{Hom}_{H_{k}%
}\left(  R_{\mu^{}_{k}}:R_{\lambda_{k}}|_{H_{k}}\right)  $, where
$\mathcal{H}_{\mu^{}_{k}}$ (resp.\ $\mathcal{H}_{\lambda_{k}}$) denotes the
representation space of $R_{\mu^{}_{k}}$ (resp.\ $R_{\lambda_{k}}$), then
since $\mathcal{H}_{\mu^{}_{k}}\subset\mathcal{I}_{n\times k}^{\left(
\mu\right)  _{k}}$ and $\mathcal{H}_{\lambda_{k}}\subset\mathcal{I}_{n\times
k}^{\left(  \lambda\right)  _{k}}$ it follows that we have an inductive chain
of homomorphisms $\left\{  T_{k}\colon\mathcal{H}_{\mu^{}_{k}}\rightarrow
\mathcal{H}_{\lambda_{k}}\right\}  $. Let $\mathcal{H}_{\mu^{}_{\infty}}$
(resp.\ $\mathcal{H}_{\lambda_{\infty}}$) denote the inductive limit of
$\mathcal{H}_{\mu^{}_{k}}$ (resp.\ $\mathcal{H}_{\lambda_{k}}$); then there
exists a unique homomorphism $T_{\infty}\colon\mathcal{H}_{\mu^{}_{\infty}%
}\rightarrow\mathcal{H}_{\lambda_{\infty}}$ (see, for example, \cite[Theorem
2.5, p.~430]{Dug78}, or \cite[p.~44]{Rot79}). Again by Theorem \ref{ThmFID.1},
$R_{\lambda_{\infty}}=\varinjlim R_{\lambda_{k}}$ (resp.\ $R_{\mu^{}_{\infty}%
}=\varinjlim R_{\mu^{}_{k}}$) is irreducible with signature $\left(
\lambda\right)  _{\infty}$ (resp.\ $\left(  \mu\right)  _{\infty}$), and it is
easy to show that $T_{\infty}$ is an intertwining homomorphism. Conversely,
all homomorphisms of inductive limits arise that way. Consequently, the chain
$\operatorname*{Hom}_{H_{k}}\left(  R_{\mu^{}_{k}}:R_{\lambda_{k}}|_{H_{k}%
}\right)  $ induces the inductive limit $\operatorname*{Hom}_{H_{\infty}%
}\left(  R_{\mu^{}_{\infty}}:R_{\lambda_{\infty}}|_{H_{\infty}}\right)  $.
Obviously for sufficiently large $k$, $\dim\left[  \operatorname*{Hom}_{H_{k}%
}\left(  R_{\mu^{}_{k}}:R_{\lambda_{k}}|_{H_{k}}\right)  \right]  =\dim\left[
\operatorname*{Hom}_{H_{\infty}}\left(  R_{\mu^{}_{\infty}}:R_{\lambda
_{\infty}}|_{H_{\infty}}\right)  \right]  $. By duality, we obtain in the same
way the inductive limit $\left[  \operatorname*{Hom}_{G_{n}^{\prime}}\left(
R_{\left(  \lambda_{n}^{\prime}\right)  _{\infty}}^{\prime}:R_{\left(  \mu
_{n}^{\prime}\right)  _{\infty}}^{\prime}\bigg|_{G_{n}^{\prime}}\right)
\right]  $; actually this chain stabilizes for $k$ sufficiently large. It
follows from Theorem \ref{ThmFD.2} (see also the proof of Theorem 4.1 in
\cite{Ton95}) that $\dim\left[  \operatorname*{Hom}_{H_{\infty}}\left(  \mu
^{}_{H_{\infty}}:\lambda_{G_{\infty}}|_{H_{\infty}}\right)  \right]
=\dim\left[  \operatorname*{Hom}_{G_{n}^{\prime}}\left(  \lambda_{\left(
G_{n}^{\prime}\right)  _{\infty}}^{\prime}:\mu_{\left(  H_{n}^{\prime}\right)
_{\infty}}^{\prime}\bigg|_{G_{n}^{\prime}}\right)  \right]  $.
\end{proof}

As an example we again consider the case $\mathcal{F}_{1\times\infty}$ with
$G_{\infty}=\mathrm{U}\left(  \infty\right)  $, $G_{1}^{\prime}=\mathrm{U}%
\left(  1\right)  $, $H_{\infty}=\mathrm{SO}\left(  \infty\right)  $, and
$H_{1}^{\prime}=\widetilde{\mathrm{SL}_{2}\left(  \mathbb{R}\right)  }$. Then
from Eq.\ (\ref{eqFD.17}), $\left(  \lambda\right)  _{k}=(\underset
{k}{\underbrace{m,0,\dots,0}})$, $\lambda_{1}^{\prime}=\left(  m\right)  $,
and $\mathcal{I}_{1\times k}^{\left(  \lambda\right)  _{k}}=\mathcal{P}%
_{1\times k}^{\left(  m\right)  }$. It follows that $\left(  \lambda\right)
_{\infty}=\left(  m,0,0,\vec{0}\right)  $ and $\mathcal{I}_{1\times\infty
}^{\left(  \lambda\right)  _{\infty}}=\mathcal{P}_{1\times\infty}^{\left(
m\right)  _{\infty}}$, the vector space of all homogeneous polynomials of
degree $m$ in infinitely many variables $Z_{1}$, $Z_{2}$, etc. The
infinitesimal action of $R_{\left(  \mathrm{U}\left(  1\right)  \right)  _{k}%
}^{\prime}$ is given by Eq.\ (\ref{eqFD.10}), $L_{k}=\sum_{i=1}^{k}%
Z_{i}\partial/\partial Z_{i}$, so the infinitesimal action $L_{\left(
m\right)  _{\infty}}$ is given the formal series $\sum_{i=1}^{\infty}%
Z_{i}\partial/\partial Z_{i}$. For $H_{\infty}=\mathrm{SO}\left(
\infty\right)  $ and $H_{1}^{\prime}=\widetilde{\mathrm{SL}_{2}\left(
\mathbb{R}\right)  }$ the actions are more delicate to describe. From
Eq.\ (\ref{eqFD.22}), $\left(  \mu\right)  _{k}=(\underset{\left[  k/2\right]
}{\underbrace{r,0,\dots,0}})$, where $r$ is an integer $\geq0$, and therefore
$\left(  \mu\right)  _{\infty}=\left(  r,0,0,\vec{0}\right)  $. Let
$\mathcal{H}_{1\times k}^{\left(  r\right)  _{k}}$ denote the space of all
harmonic homogeneous polynomials of degree $r$ in $k$ variables $Z_{1}%
,\dots,Z_{k}$ then from Eq.\ (\ref{eqFD.22}) $\mathcal{I}_{1\times k}^{\left(
r\right)  _{k}}=\sum\limits_{j=0}^{\infty}{\!{\oplus}\,}p_{0,k}^{j}%
\mathcal{H}_{1\times k}^{\left(  r\right)  _{k}}$, where $p_{0,k}\left(
Z\right)  =\sum_{i=1}^{k}Z_{i}^{2}$. We define the actions $R_{\mathrm{SO}%
\left(  \infty\right)  }$ and $R_{\left(  \widetilde{\mathrm{SL}_{2}\left(
\mathbb{R}\right)  }\right)  _{\infty}}^{\prime}$ as follows:

Consider the algebras $\left(  \mathrm{sl}_{2}\left(  \mathbb{R}\right)
\right)  _{k}$ with the bases $\left\{  E_{k}^{{}},X_{k}^{+},X_{k}%
^{-}\right\}  $ given by Eq.\ (\ref{eqFID.5}); define the \emph{projective} or
\emph{inverse limit} of the family $\left\{  \left(  \mathrm{sl}_{2}\left(
\mathbb{R}\right)  \right)  _{k},\mathcal{I}_{1\times k}^{\left(  r\right)
_{k}}\right\}  $ as follows: For each pair of indices $l,k$ with $l\leq k$ a
continuous homomorphism $\phi_{l}^{k}\colon\left(  \mathrm{sl}_{2}\left(
\mathbb{R}\right)  \right)  _{k}\rightarrow\left(  \mathrm{sl}_{2}\left(
\mathbb{R}\right)  \right)  _{l}$ by sending $E_{k}$ to $E_{l}$, $X_{k}^{+}$
to $X_{l}^{+}$, $X_{k}^{-}$ to $X^{-l}$, and extends by linearity to $\left(
\mathrm{sl}_{2}\left(  \mathbb{R}\right)  \right)  _{k}\rightarrow\left(
\mathrm{sl}_{2}\left(  \mathbb{R}\right)  \right)  _{l}$. Clearly $\phi
_{l}^{k}$ satisfies fhe following:\renewcommand{\theenumi}{\alph{enumi}}

\begin{enumerate}
\item $\phi_{k}^{k}$ is the identity map for all $k$,

\item  if $i\leq l\leq k$ then $\phi_{i}^{k}=\phi_{i}^{l}\circ\phi_{l}^{k}$.
\end{enumerate}

\noindent The inverse limit of the system $\left\{  \mathrm{sl}_{2}\left(
\mathbb{R}\right)  _{k}\right\}  $ is denoted by%
\begin{multline}
\mathrm{sl}_{2}\left(  \mathbb{R}\right)  _{\infty}=\varprojlim\mathrm{sl}%
_{2}\left(  \mathbb{R}\right)  _{k}=\left\langle E_{\infty}^{{}},X_{\infty
}^{+},X_{\infty}^{-}\right\rangle ,\\
\text{where }E_{\infty}=\frac{1}{2}1_{\infty}+L_{\infty},\;X_{\infty}%
^{+}=\frac{1}{2}\sum_{i=1}^{\infty}Z_{i}^{2},\;X_{\infty}^{-}=\frac{1}{2}%
\sum_{i=1}^{\infty}\frac{\partial\;}{\partial Z_{i}^{2}}. \label{eqFID.7}%
\end{multline}
Then $\left\{  E_{\infty}^{{}},X_{\infty}^{+},X_{\infty}^{-}\right\}  $ acts
on $\mathcal{F}_{1\times\infty}$ as follows: If $f\in\mathcal{F}%
_{1\times\infty}$ then we may assume that $f\in\mathcal{F}_{1\times k}$ for
some $k$ and%
\begin{equation}
E_{\infty}f=E_{k}f,\qquad X_{\infty}^{+}f=X_{k}^{+}f\text{\qquad and\qquad
}X_{\infty}^{-}f=X_{k}^{-}f. \label{eqFID.8}%
\end{equation}
If $\mathcal{H}_{1\times\infty}^{\left(  r\right)  _{\infty}}$ denotes the
subspace (of $\mathcal{P}_{1\times\infty}^{\left(  r\right)  _{\infty}}$) of
all harmonic homogeneous polynomials of infinitely many variables $Z_{1}$,
$Z_{2}$, etc.\ (i.e., $h\in\mathcal{H}_{1\times\infty}$ if and only if
$h\in\mathcal{P}_{1\times\infty}^{\left(  r\right)  _{\infty}}$ and
$X_{\infty}^{-}h=0$) then%
\begin{equation}
\mathcal{I}_{1\times\infty}^{\left(  r\right)  _{\infty}}=\sum\limits_{j=0}%
^{\infty}{\!{\oplus}\,}\left(  2X_{\infty}^{+}\right)  ^{j}\mathcal{H}%
_{1\times\infty}^{\left(  r\right)  _{\infty}}, \label{eqFID.9}%
\end{equation}
where in Eq.\ (\ref{eqFID.9}) $2X_{\infty}^{+}=\left(  p_{0}\right)  _{\infty
}=\sum_{i=1}^{\infty}Z_{i}^{2}$. Note that $\mathcal{H}_{1\times\infty
}^{\left(  r\right)  _{\infty}}$ corresponds to the inductive limit of the
chain $\left\{  \mathcal{H}_{1\times k}^{\left(  r\right)  _{k}}\right\}  $.
Let $R_{\mathrm{SO}\left(  \infty\right)  }^{\left(  r\right)  _{\infty}}$
denote the inductive limit representation of the chain $R_{\mathrm{SO}\left(
k\right)  }^{\left(  r\right)  _{l}}$; then $R_{\mathrm{SO}\left(
\infty\right)  }^{\left(  r\right)  _{\infty}}$ together with
Eq.\ (\ref{eqFID.8}) describes completely the action of the dual pair $\left(
\widetilde{\mathrm{SL}_{2}\left(  \mathbb{R}\right)  },\mathrm{SO}\left(
\infty\right)  \right)  $ on the isotypic component $\mathcal{I}%
_{1\times\infty}^{\left(  r\right)  _{\infty}}$ and thus we have the isotypic
decompositions for the dual pairs $\left(  \mathrm{U}\left(  1\right)
,\mathrm{U}\left(  \infty\right)  \right)  $ and $\left(  \widetilde
{\mathrm{SL}_{2}\left(  \mathbb{R}\right)  },\mathrm{SO}\left(  \infty\right)
\right)  $,%
\[
\mathcal{F}_{1\times\infty}=\sum\limits_{\makebox[0pt]{\hss$\scriptstyle
m=0$\hss}}^{\infty}{\!{\oplus}\,}\mathcal{I}_{1\times\infty}^{\left(
m\right)  _{\infty}}= \sum\limits_{\makebox[0pt]{\hss$\scriptstyle
r=0$\hss}}^{\infty}{\!{\oplus}\,}\mathcal{I}_{1\times\infty}^{\left(
r\right)  _{\infty}},
\]
and thus Theorem \ref{ThmFID.2} is verified for this example.

Since the next two examples are very important by their applications to
Physics we shall state them as corollaries to Theorem \ref{ThmFID.2}.

\begin{corollary}
\label{CorFID.3}Let $G_{\infty}$ denote the direct product of $r$ copies of
$H_{\infty}$ where $H_{\infty}=\mathrm{U}\left(  \infty\right)  $,
$\mathrm{SO}\left(  \infty\right)  $, or $\mathrm{Sp}\left(  \infty\right)  $.
If $G_{\infty}$ acts as the exterior tensor product representation $V^{\left(
\lambda_{1}\right)  _{\infty}}\otimes\dots\otimes V^{\left(  \lambda
_{r}\right)  _{\infty}}$, where each $V^{\left(  \lambda_{i}\right)  _{\infty
}}$, $1\leq i\leq r$, is an irreducible unitary $H_{\infty}$-module, then
$H_{\infty}$ acts as the inner \textup{(}or Kronecker\textup{)} tensor
product representation on $V^{\left(  \lambda_{1}\right)  _{\infty}}%
\mathop{\hat{\otimes}}\dotsb\mathop{\hat{\otimes}}V^{\left(  \lambda
_{r}\right)  _{\infty}}$. If $\lambda_{G_{\infty}}$ denotes an irreducible
unitary representation of class $\left(  \lambda_{1}\right)  _{G_{\infty}%
}\otimes\dots\otimes\left(  \lambda_{r}\right)  _{G_{\infty}}$ and $\mu
^{}_{H_{\infty}}$ denotes an irreducible unitary representation of class
$\left(  \mu\right)  _{H_{\infty}}$ then the multiplicity $\dim\left[
\operatorname*{Hom}_{H_{\infty}}\left(  \mu^{}_{H_{\infty}}:\lambda
_{G_{\infty}}|_{H_{\infty}}\right)  \right]  $ of the representation $\left(
\mu\right)  _{H_{\infty}}$ in the inner tensor product $\left(  \lambda
_{1}\right)  _{\infty}\mathop{\hat{\otimes}}\dotsb\mathop{\hat{\otimes}%
}\left(  \lambda_{r}\right)  _{\infty}$ is equal to the multiplicity of
$\left(  \mu\right)  _{H_{k}}$ in the inner tensor product $\left(
\lambda_{1}\right)  _{k}\mathop{\hat{\otimes}}\dotsb\mathop{\hat{\otimes}%
}\left(  \lambda_{r}\right)  _{k}$ for sufficiently large $k$.
\end{corollary}

\begin{proof}
If $\left(  \lambda_{i}\right)  _{\infty}=\left(  \lambda_{i}^{1},\lambda
_{i}^{2},\dots,\lambda_{i}^{j},\dots\right)  $ where $\lambda_{i}^{j}$ are
integers such that $\lambda_{i}^{1}\geq\lambda_{i}^{2}\geq\cdots$ and
$\lambda_{i}^{j}=0$ for all but a finite number of $j$, let $n$ denote the
total number of all nonzero entries $\lambda_{i}^{j}$, $1\leq i\leq r$; then
$V^{\left(  \lambda_{1}\right)  _{\infty}}\otimes\dots\otimes V^{\left(
\lambda_{r}\right)  _{\infty}}$ can be realized as a subspace of the
Bargmann--Segal--Fock space $\mathcal{F}_{n\times\infty}$. From Theorem
\ref{ThmFID.2} it follows that $V^{\left(  \lambda_{1}\right)  _{\infty}%
}\otimes\dots\otimes V^{\left(  \lambda_{r}\right)  _{\infty}}$ belongs to the
isotypic component $\mathcal{I}_{n\times\infty}^{\left(  \lambda\right)
_{G_{\infty}}}$ of $\mathcal{F}_{n\times\infty}$, thus $V^{\left(  \lambda
_{1}\right)  _{\infty}}\otimes\dots\otimes V^{\left(  \lambda_{r}\right)
_{\infty}}$ is the inductive limit of the chain $\left\{  V^{\left(
\lambda_{1}\right)  _{k}}\otimes\dots\otimes V^{\left(  \lambda_{r}\right)
_{k}}\right\}  $. If $\mu^{}_{H_{\infty}}$ is an irreducible unitary
representation of class $\left(  \mu\right)  _{H_{\infty}}$ then by Theorem
\ref{ThmFID.2}%
\[
\dim\left[  \operatorname*{Hom}\nolimits_{H_{\infty}}\left(  \mu^{}%
_{H_{\infty}}:\lambda_{G_{\infty}}|_{H_{\infty}}\right)  \right]  =\dim\left[
\operatorname*{Hom}\nolimits_{G_{n}^{\prime}}\left(  \lambda_{\left(
G_{n}^{\prime}\right)  _{\infty}}^{\prime}:\mu_{\left(  H_{n}^{\prime}\right)
_{\infty}}^{\prime}\bigg|_{G_{n}^{\prime}}\right)  \right]  ,
\]
where $\lambda_{\left(  G_{n}^{\prime}\right)  _{\infty}}^{\prime}$
(resp.\ $\mu_{\left(  H_{n}^{\prime}\right)  _{\infty}}^{\prime}$) is the
representation of $G_{n}^{\prime}$ (resp.\ $H_{n}^{{}}$) dual to
$\lambda_{G_{\infty}}$(resp.\ $\mu^{}_{H_{\infty}}$). For sufficiently large
$k$ every $\mu^{}_{H_{\infty}}$ is the inductive limit of a chain $\mu
^{}_{H_{k}}$ and for such a $k$ Theorem \ref{ThmFD.2} implies that%
\begin{align*}
\dim\left[  \operatorname*{Hom}\nolimits_{H_{k}}\left(  \mu^{}_{H_{k}}%
:\lambda_{G_{k}}|_{H_{k}}\right)  \right]   &  =\dim\left[
\operatorname*{Hom}\nolimits_{G_{n}^{\prime}}\left(  \lambda_{\left(
G_{n}^{\prime}\right)  _{k}}^{\prime}:\mu_{\left(  H_{n}^{\prime}\right)
_{k}}^{\prime}\bigg|_{G_{n}^{\prime}}\right)  \right] \\
&  =\dim\left[  \operatorname*{Hom}\nolimits_{G_{n}^{\prime}}\left(
\lambda_{\left(  G_{n}^{\prime}\right)  _{\infty}}^{\prime}:\mu_{\left(
H_{n}^{\prime}\right)  _{\infty}}^{\prime}\bigg|_{G_{n}^{\prime}}\right)
\right]  ,
\end{align*}
and this achieves the proof of Corollary \ref{CorFID.3}.
\end{proof}

\begin{remark}
The reason that this corollary only holds for sufficiently large $k$ can be
seen in the following example. Let $G_{k}=\underset{4\text{ times}%
}{\underbrace{\mathrm{U}\left(  k\right)  \times\dots\times\mathrm{U}\left(
k\right)  }}$ and $H_{k}=\mathrm{U}\left(  k\right)  $ and consider the tensor
product $(\underset{k}{\underbrace{1,0,\dots,0}})\otimes(\underset
{k}{\underbrace{2,0,\dots,0}})\otimes(\underset{k}{\underbrace{2,0,\dots,0}%
})\otimes(\underset{k}{\underbrace{3,0,\dots,0}})$; then for $k=2$ we have the
spectral decomposition
\[
\left(  1,0\right)  \otimes\left(  2,0\right)  \otimes\left(  2,0\right)
\otimes\left(  3,0\right)  =\left(  8,0\right)  +3\left(  7,1\right)
+5\left(  6,2\right)  +5\left(  5,3\right)  +2\left(  4,4\right)  ,
\]
for $k=3$ we have
\begin{multline*}
\left(  1,0,0\right)  \otimes\left(  2,0,0\right)  \otimes\left(
2,0,0\right)  \otimes\left(  3,0,0\right) \\
=\left(  8,0,0\right)  +3\left(  7,1,0\right)  +5\left(  6,2,0\right)
+5\left(  5,3,0\right)  +2\left(  4,4,0\right) \\
+3\left(  6,1,1\right)  +6\left(  5,2,1\right)  +5\left(  4,3,1\right)
+3\left(  4,2,2\right)  +2\left(  3,3,2\right)  ,
\end{multline*}
for $k\geq4$ we have
\begin{multline*}
(\underset{k}{\underbrace{1,0,\dots,0}})\otimes(\underset{k}{\underbrace
{2,0,\dots,0}})\otimes(\underset{k}{\underbrace{2,0,\dots,0}})\otimes
(\underset{k}{\underbrace{3,0,\dots,0}})\\
=\left(  8,0,\dots,0\right)  +3\left(  7,1,0,\dots,0\right)  +5\left(
6,2,0,\dots,0\right)  +5\left(  5,3,0,\dots,0\right) \\
+2\left(  4,4,0,\dots,0\right)  +3\left(  6,1,1,0,\dots,0\right)  +6\left(
5,2,1,0,\dots,0\right)  +5\left(  4,3,1,0,\dots,0\right) \\
+3\left(  4,2,2,0,\dots,0\right)  +2\left(  3,3,2,0,\dots,0\right)  +\left(
5,1,1,1,0,\dots,0\right) \\
+2\left(  4,2,1,1,0,\dots,0\right)  +\left(  3,3,1,1,0,\dots,0\right)
+\left(  3,2,2,1,0,\dots,0\right)  .
\end{multline*}
Thus we can see that the spectral decomposition of $\left(  1,\vec{0}\right)
_{\infty}\otimes\left(  2,\vec{0}\right)  _{\infty}\otimes\left(  2,\vec
{0}\right)  _{\infty}\otimes\left(  3,\vec{0}\right)  _{\infty}$ is the same
as that of order $k$ for $k\geq4$, with infinitely many zeroes at the end of
each signature.

Note also that this corollary applied to the tensor product $\underset{r\text{
times}}{\underbrace{\left(  1,\vec{0}\right)  _{\infty}\otimes\dots
\otimes\left(  1,\vec{0}\right)  _{\infty}}}$ together with the Schur--Weyl
Duality Theorem for $\mathrm{U}\left(  r\right)  $ implies the generalized
Schur--Weyl Duality Theorem proved by Kirillov for $\mathrm{U}\left(
\infty\right)  $ in \cite{Kir73}.
\end{remark}

\begin{corollary}
\label{CorFID.4}Let $V^{\left(  \lambda_{1}\right)  _{\infty}},\dots
,V^{\left(  \lambda_{r}\right)  _{\infty}}$ and $V^{\left(  \mu\right)
_{\infty}}$ be irreducible unitary representation of $H_{\infty}$. Let
$V^{\left(  \mu^{\checkmark}\right)  _{\infty}}$ be the representation
\textup{(}of $H_{\infty}$\textup{)} contragredient to $V^{\left(  \mu\right)
_{\infty}}$. Let $I^{\infty}$ denote the equivalence class of the identity
representation of $H_{\infty}$. Then the multiplicity of $\left(  \mu\right)
_{\infty}$ in the tensor product $\left(  \lambda_{1}\right)  _{\infty}%
\mathop{\hat{\otimes}}\dotsb\mathop{\hat{\otimes}}\left(  \lambda_{r}\right)
_{\infty}$ is equal to the multiplicity of $I^{\infty}$ in the tensor product
$\left(  \lambda_{1}\right)  _{\infty}\mathop{\hat{\otimes}}\dotsb\mathop
{\hat{\otimes}}\left(  \lambda_{r}\right)  _{\infty}\mathop{\hat{\otimes}%
}\left(  \mu^{\checkmark}\right)  _{\infty}$.
\end{corollary}

\begin{proof}
To prove this corollary we apply Corollary \ref{CorFID.3} to $G_{\infty
}=\underset{r}{\underbrace{H_{\infty}\times\cdots\times H_{\infty}}}$ and
$G_{k}=\underset{r}{\underbrace{H_{k}\times\cdots\times H_{k}}}$, then apply
Theorem \ref{ThmFID.2} to $G_{\infty}=\underset{r}{\underbrace{H_{\infty
}\times\cdots\times H_{\infty}}}\times H_{\infty}^{\checkmark}$ and
$G_{k}=\underset{r}{\underbrace{H_{k}\times\cdots\times H_{k}}}\otimes
H_{k}^{\checkmark}$, and finally apply Theorem 2.1 of \cite{KlTo96} to obtain
the desired result at order $k$. The main difficulty resides with the
definition of the identity representation on $V^{\left(  \lambda_{1}\right)
_{\infty}}\mathop{\hat{\otimes}}\dotsb\mathop{\hat{\otimes}}V^{\left(
\lambda_{r}\right)  _{\infty}}\mathop{\hat{\otimes}}V^{\left(  \mu
^{\checkmark}\right)  _{\infty}}$, which we will construct below.

For each $k$ let $I^{k}$ denote the identity representation of $H_{k}$ on
\linebreak $V^{\left(  \lambda_{1}\right)  _{k}}\mathop{\hat{\otimes}}%
\dotsb\mathop
{\hat{\otimes}}V^{\left(  \lambda_{r}\right)  _{k}}\mathop{\hat{\otimes}%
}V^{\left(  \mu^{\checkmark}\right)  _{k}}$. This means that if $I^{k}$ occurs
with multiplicity $d$ in $V^{\left(  \lambda_{1}\right)  _{k}}\mathop
{\hat{\otimes}}\dotsb\mathop{\hat{\otimes}}V^{\left(  \lambda_{r}\right)
_{k}}\mathop{\hat{\otimes}}V^{\left(  \mu^{\checkmark}\right)  _{k}}$ then
there exist $d$ nonzero vectors $f_{i,k}$, $i=1,\dots,d$, such that $R_{H_{k}%
}\left(  u\right)  f_{i,k}=f_{i,k}$ for all $u\in H_{k}$. By construction each
$f_{i,k}$ is a polynomial function in $\mathcal{F}_{n\times k}$ for some $n$.
Thus $f_{i,k}$ is an $H_{k}$-invariant polynomial in $\mathcal{F}_{n\times k}%
$. If $J_{i,k}$ denotes the one-dimensional subspace spanned by $f_{i,k}$,
then for sufficiently large $k$ and for each fixed $i=1,\dots,d$ we have a
chain of irreducible unitary representations $\left\{  H_{k},I^{k}%
,J_{i,k}\right\}  _{k}$. We can define the isomorphism $\psi_{k+1}^{k}\colon
J_{i,k}\rightarrow J_{i,k+1}$ by $\psi_{k+1}^{k}\left(  cf_{i,k}\right)
=cf_{i,k+1}$, $c\in\mathbb{C}$; then obviously
\[
\psi_{k+1}^{k}\left(  R_{H_{k}}\left(  u\right)  f_{i,k}\right)  =R_{H_{k+1}%
}\left(  u\right)  f_{i,k+1}=R_{H_{k+1}}\left(  u\right)  \psi_{k+1}%
^{k}\left(  f_{i,k}\right)  ,
\]
for all $u\in H_{k}$. Also for all $k$, $l$, $m$ with $k\leq l\leq m$ we have
$\psi_{m}^{k}=\psi_{m}^{l}\circ\psi_{l}^{k}$. Thus we can define the inductive
limit representation $\left\{  H_{\infty},I^{\infty},J_{i,\infty}\right\}  $,
where the action of $H_{\infty}$ on $J_{i,\infty}$ is defined as follows:

Let $u\in H_{\infty}$; then $u\in H_{k}$ for some $k$. If $f\in J_{i,l}$ for
some $l$ then%
\[
R_{H_{\infty}}\left(  u\right)  f_{l}=R_{H_{k}}\left(  u\right)  \psi_{k}%
^{l}f\text{\qquad for }l<k,
\]
and%
\[
R_{H_{\infty}}\left(  u\right)  f_{l}=R_{H_{k}}\left(  u\right)  \psi_{l}%
^{k}f\text{\qquad for }k\leq l.
\]
Then it follows from Theorem \ref{ThmFID.1} that $\left\{  H_{\infty
},I^{\infty},J_{i,\infty}\right\}  $ is irreducible with signature $\left(
\vec{0}\right)  _{\infty}$. The only problem with this approach is that the
isomorphism embedding $\psi_{k+1}^{k}$ is not the isomorphic embedding
$i_{k+1}^{k}\colon\mathcal{F}_{n\times k}\rightarrow\mathcal{F}_{n\times
\left(  k+1\right)  }$. To circumvent this difficulty we define the
\emph{inverse} or \emph{projective limit} of the family $\left\{  H_{k}%
,I^{k},J_{k}\right\}  $ where $J_{k}$ denotes the subspace of all $H_{k}%
$-invariants in $V^{\left(  \lambda_{1}\right)  _{k}}\mathop{\hat{\otimes}%
}\dotsb\mathop{\hat{\otimes}}V^{\left(  \lambda_{r}\right)  _{k}}\mathop
{\hat{\otimes}}V^{\left(  \mu^{\checkmark}\right)  _{k}}$, as follows: For
each pair of indices $l$, $k$ with $l\leq k$ define a continuous homomorphism
$\phi_{l}^{k}\colon J_{k}\rightarrow J_{l}$ such that\renewcommand{\theenumi
}{\roman{enumi}}

\begin{enumerate}
\item \label{CorFID.4proof(1)}$\phi_{k}^{k}$ is the identity map on $J_{k}$,

\item \label{CorFID.4proof(2)}if $i\leq l\leq k$ then $\phi_{i}^{k}=\phi
_{i}^{l}\circ\phi_{l}^{k}$.
\end{enumerate}

\noindent Here we can take $\phi_{l}^{k}$ as the \emph{truncation}
homomorphism, i.e., $\phi_{l}^{k}$ is defined on the generators $f_{i,k}$ by%
\[
\phi_{l}^{k}\left(  f_{i,k}\right)  =f_{i,l}.
\]
The \emph{projective limit} of the system $\left\{  H_{k},J_{k},\phi_{l}%
^{k}\right\}  $ is then formally defined by%
\[
J_{\infty_{\leftarrow}}:=\varprojlim J_{k}=\left\{  \left(  f_{k}\right)
\in\prod_{k}J_{k}:f_{l}=\phi_{l}^{k}\left(  f_{k}\right)  ,\;\forall\,l\leq
k\right\}  .
\]
Let $\pi_{k}\colon J_{\infty_{\leftarrow}}\rightarrow J_{k}$ denote the
projection of $J_{\infty_{\leftarrow}}$ onto $J_{k}$. Let $I^{\infty
_{\leftarrow}}$ denote the representation of $H_{\infty}$ on $J_{\infty
_{\leftarrow}}$; then $\pi_{k}\left(  I^{\infty_{\leftarrow}}f\right)
=\pi_{k}\left(  f\right)  $. Recall that if $\mathcal{P}_{n\times k}$ denotes
the subspace of all polynomial functions on $\mathbb{C}^{n\times k}$ the
$\mathcal{P}_{n\times k}$ is dense in $\mathcal{F}_{n\times k}$. Let
$\mathcal{P}_{n\times\infty}=\bigcup_{k=1}^{\infty}\mathcal{P}_{n\times k}$
denote the inductive limit of $\mathcal{P}_{n\times k}$; then clearly
$\mathcal{P}_{n\times\infty}$ is dense in $\mathcal{F}_{n\times\infty}$. Let
$\mathcal{P}_{n\times\infty}^{\ast}$ (resp.\ $\mathcal{F}_{n\times\infty
}^{\ast}$) denote the \emph{dual} or \emph{adjoint} space of $\mathcal{P}%
_{n\times\infty}$ (resp.\ $\mathcal{F}_{n\times\infty}$). Then since
$\mathcal{P}_{n\times\infty}$ is dense in $\mathcal{F}_{n\times\infty}$,
$\mathcal{F}_{n\times k}^{\ast}$ is dense in $\mathcal{P}_{n\times\infty
}^{\ast}$. By the Riesz representation theorem for Hilbert spaces, every
element $f^{\ast}\in\mathcal{F}_{n\times\infty}^{\ast}$ is of the form
$\ip{\,\cdot\,}{f}$ for some $f\in\mathcal{F}_{n\times\infty}$, and the map
$f^{\ast}\rightarrow f$ is an anti-linear (or conjugate-linear) isomorphism.
Thus we can identify $\mathcal{F}_{\infty}^{\ast}$ with $\mathcal{F}_{\infty
}^{{}}$ and obtain the \emph{rigged Hilbert space} as the triple
$\mathcal{P}_{n\times\infty}^{{}}\subset\mathcal{F}_{n\times\infty}^{{}%
}\subset\mathcal{P}_{n\times\infty}^{\ast}$ (see \cite{GeVi64} for the
definition of rigged Hilbert spaces). However, generally an element of
$J_{\infty_{\leftarrow}}$ does not belong to $\mathcal{P}_{n\times\infty
}^{\ast}$, but can still be considered as a linear functional (not necessarily
continuous) on $\mathcal{P}_{n\times\infty}$, and furthermore, in this context
the identity representation $I^{\infty_{\leftarrow}}$ will respect the
isomorphic embedding $i_{k+1}^{k}\colon\mathcal{F}_{n\times k}\rightarrow
\mathcal{F}_{n\times\left(  k+1\right)  }$.
\end{proof}

\section{\label{Con}Conclusion}

We have studied thoroughly several reciprocity theorems for some dual pairs of
groups $\left(  G_{n}^{\prime},G_{\infty}^{{}}\right)  $ and $\left(
H_{n}^{\prime},H_{\infty}^{{}}\right)  $, where $G_{\infty}$ is the inductive
limit of a chain $\left\{  G_{k}\right\}  $ of compact groups, $H_{\infty}$ is
the inductive limit of a chain $\left\{  H_{k}\right\}  $ such that for each
$k$, $H_{k}$ is a compact subgroup of $G_{k}$, and $G_{n}^{\prime}\subset
H_{n}^{\prime}$ are finite-dimensional Lie groups. These theorems show, in
particular, that the multiplicity of an irreducible unitary representation of
$H_{\infty}$ with signature $\left(  \mu\right)  _{H_{\infty}}$ in the
restriction to $H_{\infty}$ of an irreducible unitary representation of
$G_{\infty}$ with signature $\left(  \lambda\right)  _{G_{\infty}}$ is always
finite. This is extremely important in the problem of spectral decompositions
of tensor products of irreducible unitary representations of inductive limits
of compact classical groups. This type of problems arises naturally in Physics
(cf.\ \cite{KaRa87}), and in \cite{HoTo98} tensor product decompositions of
tame representations of $\mathrm{U}\left(  \infty\right)  $ are investigated.
In \cite{Ol'90} Ol'shanskii generalized Howe's theory of dual pairs to some
infinite-dimensional dual pairs of groups. This is the right context to
generalize the reciprocity theorem \ref{ThmFID.2} for the infinite-dimensional
dual pairs $\left(  G_{\infty}^{\prime},G_{\infty}^{{}}\right)  $ and $\left(
H_{\infty}^{\prime},H_{\infty}^{{}}\right)  $ which will be part of our work
in a forthcoming publication.

%\bibliographystyle{BFTALPHA}
%\bibliography{tonthat}
%
%
%
%
%
%
%

\end{document}